\documentclass[12pt]{amsart}

\usepackage{amsthm,amsfonts, amssymb, amscd, mathtools}
\usepackage{youngtab}
\usepackage{ytableau}
\usepackage[all]{xy}
\usepackage{euscript}
\usepackage{cite}
\usepackage{comment}

\usepackage{tikz}
\usetikzlibrary{matrix,arrows, decorations.markings}
\usetikzlibrary{positioning}

\usepackage[normalem]{ulem}

\usepackage{hyperref}
\hypersetup{%
    linktoc=page
}

\let\oldtocsection=\tocsection
\let\oldtocsubsection=\tocsubsection
\renewcommand{\tocsection}[2]{\hspace{0em}\oldtocsection{#1}{#2}}
\renewcommand{\tocsubsection}[2]{\hspace{1em}\oldtocsubsection{#1}{#2}}
\newtheorem{theorem}{Theorem}[section]
\newtheorem{lemma}[theorem]{Lemma}
\newtheorem{corollary}[theorem]{Corollary}
\newtheorem{proposition}[theorem]{Proposition}

\newtheorem{conjecture}[theorem]{Conjecture}
\theoremstyle{remark}
\newtheorem{remark}[theorem]{Remark}
\theoremstyle{remark}
\newtheorem{example}[theorem]{Example}
\theoremstyle{definition}

\theoremstyle{definition}

\theoremstyle{definition}

\newtheorem{definition}[theorem]{Definition}

\numberwithin{equation}{section}

\newcommand{\C}{\mathbb{C}}           % Use for complex numbers.
\newcommand{\N}{\mathbb{ N}}           % Use for positive integers.
\newcommand{\Z}{\mathbb{ Z}}           % Use for integers.
 
\newcommand{\cd}{\mathcal{D}}

 \newcommand{\cp}{\mathcal{P}}

\renewcommand{\tilde}{\widetilde}

\newcommand{\qchoos}[2]{\genfrac[]{0pt}{0}{#1}{#2}_q}

\newcommand{\inv}{\mathrm{inv}}

\newcommand{\hfam}{\mathbf{h}}
\newcommand{\Hinv}{I_\hfam}
\newcommand{\hinv}{\mathcal{I}_\hfam}

\newcommand{\Ext}{\mathrm{Ext}}
\newcommand{\Symm}{S}

\newcommand{\bt}{\mathbf{t}}
\newcommand{\bj}{\mathbf{j}}
\newcommand{\bm}{\mathbf{m}}

\begin{document}

\title[Inversion polynomials]{Restricted inversion polynomials}

\author{Jeongwon Lee} \address{Department of Mathematics\\ Washington University in St. Louis \\ One Brookings Drive \\ St.~Louis, Missouri  63130 \\ USA }
\email{jeongwon.l@wustl.edu}

\author{Nathan Lesnevich} \address{Department of Mathematics\\ Oklahoma State University \\ 401 Mathematical Sciences \\ Stillwater, OK 74075 \\ USA }
\email{nlesnev@okstate.edu}
\urladdr{\url{https://nlesnevich.github.io/}}

\author{Martha Precup} \address{Department of Mathematics\\ Washington University in St. Louis \\ One Brookings Drive \\ St.~Louis, Missouri  63130 \\ USA }
\email{martha.precup@wustl.edu}
\urladdr{\url{https://www.math.wustl.edu/~precup/}}

\begin{abstract} 
For a finite subset $I$ of positive integers, the descent polynomial $\cd(I;n)$ counts the number of permutations in $\Symm_n$ that have descent set $I$. We generalize descent polynomials by considering permutations with a specific subset $S$ of common inversions called $\hfam$-inversions, where $\hfam = (\hfam(1), \hfam(2), \ldots )$ is a weakly increasing sequence of positive integers such that $\hfam(i)> i$. We prove that this more general count, denoted by $\hinv(S;n)$, is also a polynomial. We give three explicit expansions for $\hinv(S;n)$, prove the coefficients for two of these expansions are log-concave, and define a graded generalization.
\end{abstract}

\maketitle

\section{Introduction} 
Let $I$ be a subset of $[n-1]$. MacMahon~\cite{MacMahon} proved that the number of permutations $\cd(I;n)$ in the symmetric group $\Symm_n$ with descent set equal to $I$ is a polynomial in $n$. These {descent polynomials} were subsequently studied by Diaz-Lopez, Harris, Insko, Omar, and Sagan who also gave several explicit formulas, including positive expansions in particular binomial bases~\cite{descent_poly}. Further expansions and an analysis of the resulting coefficients were given by Bencs \cite{Bencs}, and a graded generalization was studied by Gaetz and Gao \cite{Gaetz-Gao2021}.

In this paper, we generalize descent polynomials as follows. Consider a sequence $\hfam = (\hfam(1), \hfam(2), \ldots )$ of positive integers such that $\hfam(i) > i$ for all $i$. The \emph{$\hfam$-inversions} of a permutation $\pi \in \Symm_n$ is the set 
\[\inv_\hfam(\pi) \coloneqq \left\{ (i,j) \mid i < j \leq \hfam(i)\; \text{ and }\; \pi_i > \pi_j   \right\}.\]
If $\hfam(i) = i+1$ for all $i$, then the $\hfam$-inversions of $\pi$ are exactly the pairs $(i,i+1)$ such that $i$ is a descent of $\pi$. Given a subset  $S \subset \{(i,j) \mid i<j\leq \hfam(i)\}$ we consider the count:
\[
\hinv(S;n) = \#\{\pi \in \Symm_n \mid \inv_\hfam(\pi) = S\}.
\]
We say that $S$ is \emph{$\hfam$-admissible} when $\hinv(S;n)\neq 0$ for some $n$. When $\hfam$ is the sequence $\hfam(i)=i+1$ for all $i$, then $\hinv(S;n)$ is a descent polynomial.

We prove a variety of enumerative results for $\hinv(S;n)$, generalizing several known results for descent polynomials in the process. Our first main result is:
\begin{theorem}\label{intthm:poly}
Let $\hfam = (\hfam(1), \hfam(2), \ldots)$ be a sequence of positive integers such that $\hfam(i)>i$ and $S\subset \{(i,j) \mid i<j\leq \hfam(i)\}$ be an $\hfam$-admissible set. Let $\bj(S) := \max\{j \mid (i,j) \in S\}$. Then $\hinv(S;n)$ is given by a polynomial for all $n \geq \bj(S)$.
\end{theorem}
Theorem \ref{intthm:poly} is a consequence of Proposition~\ref{prop.polynomial} below. We call $\hinv(S;n)$ the \emph{$\hfam$-inversion polynomial} of $S$, or if $\hfam$ is fixed, the \emph{restricted inversion polynomial}.

In light of Theorem \ref{intthm:poly}, it is natural to ask for explicit expansions of $\hinv(S;n)$.  Explicit expansions for the descent polynomial were given by Diaz-Lopez, Harris, Insko, Omar, and Sagan\cite[Theorem 3.3]{descent_poly} and also by Bencs~\cite[Proposition 3.1]{Bencs}. 
The following generalizes an expansion of Bencs~\cite[Proposition 3.1]{Bencs} to the setting of restricted inversion polynomials (see Corollary~\ref{cor.polynomial} below).

\begin{theorem}\label{intthm:t_expansion}
Let $\hfam = (\hfam(1), \hfam(2), \ldots)$ be a sequence of positive integers such that $\hfam(i)>i$ and $S$ be an $\hfam$-admissible subset. Then there exists a sequence of non-negative integers $(t_k(S))_{k=1}^{\bj(S)}$ such that 
\begin{align*}
        \hinv(S;n) &=  \sum_{k=1}^{\bj(S)} t_k(S) {n-k \choose \bj(S)-k} 
    \end{align*}
where the values $t_k(S)$ are obtained by enumerating certain permutations with $\hfam$-inversion set equal to $S$.
\end{theorem}

When $\hfam$ is a weakly increasing sequence, also called a \emph{Hessenberg sequence}, we obtain generalizations of both previously studied expansion formulas (see Corollary~\ref{cor.k-expansion} and Theorem \ref{thm.a-expansion} below).

\begin{theorem}\label{intthm:expansion}
Let $\hfam = (\hfam(1), \hfam(2), \ldots)$ be a weakly increasing sequence of positive integers such that $\hfam(i)>i$ and $S$ be an $\hfam$-admissible subset. Let $m = \bm(S) := \max\{i \mid (i,i+1)\in S\}$ be the maximum descent in $S$. Then there exists sequences of non-negative integers $(a_k(S))_{k=0}^m$ and $(b_k(S))_{k=\hfam(m)-m}^{\hfam(m)}$ such that 
\begin{align*}
        \hinv(S;n) &=  \sum_{k=0}^m a_k(S) {n-\hfam(m)+1\choose k}\\
        &= \sum_{\ell =0}^{m} b_{\hfam(m)-\ell}(S) {n-\hfam(m)+\ell \choose \ell},
    \end{align*}
where the values $a_k(S)$ and $b_k(S)$ are each obtained by enumerating certain permutations with $\hfam$-inversion set equal to $S$.
\end{theorem}

Explicit formulas for the coefficient sequences are given in Sections~\ref{sec.polynomality} and~\ref{sec.expansion}.  Once we have identified sequences of coefficients for the restricted inversion polynomial, it is furthermore natural to ask what properties these sequences satisfy. For the descent polynomial, \cite{descent_poly} conjectured and \cite{Bencs} proved that these coefficients are log-concave. We prove the following generalizations (see Propositions~\ref{prop:t_logconcave}, \ref{prop:b_logconcave}, and~\ref{prop:a_logconcave} below).

\begin{theorem}\label{intthm:log} For every non-empty $\hfam$-admissible set $S$, the sequences $(t_k(S))$, $(a_k(S))$, and $(b_k(S))$ of coefficients from Theorems~\ref{intthm:t_expansion} and~\ref{intthm:expansion} are log-concave.
\end{theorem}

Defined by Gaetz and Gao~\cite{Gaetz-Gao2021}, the graded descent polynomial $\cd(I,n;q)$ is a polynomial in $q$ for each positive integer $n$.  Evaluation of $\cd(I,n;q)$ at $n$ gives the length-generating function in $q$ for the set of permutations with descent set $I$. 
Inspired by that work, we define a graded variation $\hinv(S,n;q)$ of the polynomial $\hinv(S;n)$ using the length statistic on permutations and compute two expansions for this graded version (Theorems~\ref{thm:t_graded} and~\ref{thm:b_graded} below). We do not fully generalize all results of Gaetz and Gao from the graded descent to graded $\hfam$-inversion polynomial. In particular, Gaetz and Gao prove that the polynomial coefficient sequence determined by the  expansion of the graded descent polynomial is strongly $q$-log concave. We do not know if the polynomial sequences of coefficients determined by the analogous expansions of $\hinv(S,n;q)$ are also strongly $q$-log concave, but have confirmed the statement is true in many examples (see Conjecture~\ref{conj:strongly_q_logconcave} below).

The rest of the paper is structured as follows. In Section 2 we provide the necessary background and definitions, and prove Theorems \ref{intthm:poly} and~\ref{intthm:t_expansion}. In Section 3 we consider the case where $\hfam$ is weakly increasing and describe and prove the two expansions for $\hinv(S;n)$ from Theorem~\ref{intthm:expansion}. In Section 4 we prove the log-concavity results of Theorem~\ref{intthm:log} by connecting our coefficient sequences $(t_k(S))$ and $(b_k(S))$ to the height polynomial of a particular poset. We define and study the graded inversion polynomial $\hinv(S,n;q)$ in Section 5.

\section{Inversion polynomials}
For each positive integer $n$ we let $S_n$ denote the symmetric group on $[n] = \{1,2,\ldots, n\}$. For any positive integer $k$ with $k\leq n$, we write $[k,n]:= \{k,k+1, \ldots, n\}$.  We denote an element $\pi\in S_n$ of the symmetric group using one-line notation by $\pi=\pi_1\pi_2\cdots \pi_n$. An \emph{inversion} of $\pi$ is a pair $(i,j)\in [n]\times[n]$ such that $i<j$ and $\pi_i>\pi_j$. We denote the set of inversions of $\pi$ by $\inv(\pi)$. We say that $i\in [n-1]$ is a \emph{descent of $\pi$} when $(i,i+1)\in \inv(\pi)$ and also call the pair $(i,i+1)$ a descent in the set $\inv(\pi)$ of inversions.

\subsection{Restricted inversion sets} Let $\hfam = (\hfam(1), \hfam(2), \hfam(3), \cdots)$ denote a sequence of positive integers such that $\hfam(i)>i$ for all $i$. We frequently assume $\hfam$ is weakly increasing, but that condition is not necessary for all of the results below. When $\hfam$ satisfies  the additional condition that $h(i) \leq h(i+1)$ for all $i$, we call $\hfam$ a \emph{Hessenberg sequence}. Motivation for this terminology comes from Remark~\ref{rem.Hessenberg} below.  Given permutation $\pi\in \Symm_n$ we define the set of \emph{$\hfam$-inversions} of $\pi$ to be 
\[
\inv_\hfam(\pi) := \{ (i,j) \in \inv(\pi) \mid  j \leq \hfam(i)  \}. 
\] 
Let $\cp_\hfam:= \{(i,j) \mid i<j\; \text{ and }\; j\leq \hfam(i)\}$ be the set of possible $\hfam$-inversions. To align our work with previous work on descent polynomials, we call pairs in $\cp_\hfam$ of the form $(i,i+1)$ descents. Since $\hfam(i)>i$ the set $\cp_\hfam$ contains all possible descents.  We say a subset $S\subset \cp_\hfam$ is \emph{$\hfam$-admissible} if $S$ is the set of $\hfam$-inversions of some permutation. Any $\hfam$-admissible set must be finite since we only consider finite permutations.

\begin{example}\label{ex:h-ad1} Let $\hfam(i)=i+2$ for all $i$. The permutation $\pi = 45231\in S_5$ has inversion set
\[
\inv(\pi) = \{(1,3),(1,4),(1,5), (2,3),(2,4),(2,5), (3,5), (4,5)\}.
\]
Since $\hfam(i)=i+2$, the set $\inv_\hfam(\pi)$ consists of the inversion pairs above which differ by at most two. We obtain
\[
\inv_{\hfam}(\pi) = \{(1,3),(2,3),(2,4),(3,5), (4,5)\}.
\]
Thus, $S=\{(1,3),(2,3),(2,4),(3,5), (4,5)\}$ is $\hfam$-admissible.
\end{example}

\begin{example}\label{ex.descents} When $\hfam(i)=i+1$ for all $i$ the set of $\hfam$-inversions of $\pi$ is the set of descents of $\pi$. Indeed, we have $i<j\leq \hfam(i)=i+1$ if and only if $j = i+1$ and 
\[
\inv_\hfam(\pi) = \{(i,i+1) \mid \pi_{i}> \pi_{i+1}\}
\] 
in this case.
\end{example}

We say that a subset $S\subseteq \cp_\hfam$ is \emph{$\hfam$-closed} if whenever $(i,j), (j,k)\in S$ with $i<j<k$ and $(i,k)\in \cp_\hfam$ we have $(i,k)\in S$.  The following fact is relatively well-known; see \cite[Proposition 6.1]{Sommers-Tymoczko2006} for a proof when $\hfam$ is weakly increasing. The proof in the more general case is straightforward.

\begin{lemma}\label{lemma.S.closed} Let $S\subseteq \cp_\hfam$ be a finite set. Then $S$ is $\hfam$-admissible if and only if both $S$ and its complement $S^c:= \cp_\hfam\setminus S$ are $\hfam$-closed subsets.
\end{lemma}

Since every non-identity permutation has at least one descent, every nonempty $\hfam$-admissible set $S$ must also contain at least one descent. We define 
\[
\bm(S) := \max\{i \mid (i,i+1)\in S\}.
\]
to be the maximum descent of $S$. 

We can now introduce the main topic of interest in this paper. Consider the set
\[
\Hinv(S,n):=\{\pi\in \Symm_n \mid  \inv_\hfam (\pi) =S\}
\]
of all permutations of $\Symm_n$ with the same $\hfam$-inversions. We also set
\begin{eqnarray}\label{eqn.def-inv-poly}
\hinv(S;n):=\#\{\pi\in \Symm_n \mid  \inv_\hfam (\pi) =S\}.
\end{eqnarray}
When $\hfam(i)=i+1$ for all $i$ the expression $\hinv(S;n)$ is the descent polynomial (cf.~Example~\ref{ex.descents}).

\begin{example} Let $\hfam(i) = i+2$ for all $i$ and $S = \{(1,2),(2,3),(1,3)\}$. Then $w \in \Hinv(S,n)$ for $n \geq 3$ if and only if $w(1) > w(2) > w(3)$ and $w(3) < w(2) < w(4) < w(5) < \cdots < w(n)$. These permutations are precisely of the form 
\[
w = k2134\cdots \hat{k}\cdots n
\]
for each $k \in \{3,\ldots,n\}$, where $\hat{k}$ indicates that $k$ is deleted from the sequence $34\ldots n$. We get that $\hinv(S;n) = \#\{3,\ldots,n\} =  n-2$.
\end{example}

As indicated in~\cite[Corollary 2.2]{descent_poly}, descent polynomials for nonempty subsets of positive integers have positive degree. 
Unlike the case of descent polynomials, it is possible for cardinality $\hinv(S;n)$ of restricted inversion sets to be constant for nonempty sets $S$. We characterize when this occurs in Corollary \ref{cor:constant} below.

\begin{example}\label{ex:constant_pre} Let $\hfam = (2,4,4,5,6,7,8,\cdots)$ and $S = \{(2,3)\}$. When $n=3$, there are precisely two permutations in $\Symm_3$ with  $\hfam$-inversion set equal to $S$, namely $132$ and $231$. When $n=4$ we cannot place $4$ in any of the first three positions without changing the $\hfam$-inversion set. Thus, $\Hinv(S,4) = \{1324, 2314\}$. This trend continues, and $\hinv(\{(2,3)\}; n) = 2$ for all $n\geq 3$. 
\end{example}

\begin{remark}\label{rem.Hessenberg} Restricted inversion sets arise in the study of the Betti numbers of certain subvarieties of the flag variety called regular semisimple Hessenberg varieties~\cite{MDeMari_Procesi_Shayman_Hess,Tymoczko_paving,Precup}. Suppose $n$ is fixed and let $h: [n]\to [n]$ be a weakly increasing function such that $h(i)> i$, called an indecomposable \emph{Hessenberg function}. Given a diagonalizable matrix $X$ with distinct eigenvalues the corresponding \emph{regular semisimple Hessenberg variety} $\mathrm{Hess}(X,h)$ consists of all nested sequences of subspaces $0\subset V_1\subset V_2\subset \cdots \subset \C^n$ such that $X(V_i)\subseteq V_{h(i)}$. The variety $\mathrm{Hess}(X,h)$ is irreducible and smooth. Let $\hfam$ be any sequence such that $\hfam(i)=h(i)$ when $i\leq n$ and $\hfam(i)>i$ otherwise.  In the language of restricted inversion polynomials,
\[
\mathrm{Poin}(\mathrm{Hess}(X,h);t) =  \sum_{\text{$S$ $\hfam$-admiss.}} \hinv(S;n)\,t^{2|S|}
\]
where $\mathrm{Poin}(\mathrm{Hess}(X,h);t)$ denotes the Poincar\'e polynomial of $\mathrm{Hess}(X,h)$. Furthermore, $\hinv(S;n)$ is the size of a particular weak left Bruhat interval in $S_n$~\cite{Harada-Precup2023}. 
\end{remark}

\subsection{Polynomality}\label{sec.polynomality}
We now turn our attention to showing that $\hinv(S;n)$ is a polynomial for each admissible set $S$. The key idea is to partition the set $\Hinv(S,n)$ using the flattening map on permutations.

Given a permutation $\pi\in \Symm_n$ and positive integer $k\leq n$ we let $\pi|_k$ denote the permutation of $\Symm_k$ obtained from $\pi$ by listing $1, 2, \ldots, k$ in the same relative order as the first $k$ entries of $\pi$. This operation is frequently called the \emph{flattening map}. For example, if $\pi = 641235$ then $\pi|_4 = 4312$ and $\pi|_3 = 321$.

For any $\hfam$-admissible subset $S\subset \cp_\hfam$, we define
\[
\bj(S):= \max\{j\in \Z_{>0} \mid (i,j)\in S \}.
\]
In other words, $\bj(S)$ is the largest index that appears in any pair of $S$. 

\begin{lemma}\label{lemma.flatten} Suppose $S$ is $\hfam$-admissible and let $n$ and $k$ be positive integers such that $\bj(S) \leq k \leq n$. If $\pi\in \Hinv(S,n)$, then $\pi|_{k} \in \Hinv(S, k)$. 
\end{lemma}
\begin{proof} By definition, $(\pi|_k)_i < (\pi|_k)_j$ if and only if $\pi_i<\pi_j$ for all $i$ and $j$ in $[k]$. Thus, 
\[
\inv_\hfam(\pi|_k) = \inv_\hfam(\pi) \cap \{(i, j) \in \cp_{\hfam}\mid i, j \in [k] \} = S \cap \{(i, j) \in \cp_{\hfam}\mid i, j \in [k] \}.
\]
Since $k\geq \bj(S)$ the intersection on the right hand side of the equation above is equal to $S$. Thus $\pi|_k \in \Hinv(S,k)$.
\end{proof}

\begin{corollary}\label{cor.smallest} The value $\bj(S)$ is the smallest positive integer $n$ such that $\hinv(S;n)> 0$.
\end{corollary}
\begin{proof}
It's obvious that $\Hinv(S,n)=\varnothing$ whenever $n<\bj(S)$ so $\hinv(S;n)=0$ for all $n<\bj(S)$. On the other hand, if $S$ is $\hfam$-admissible, then there exists some permutation $\pi\in \Symm_n$ with $n\geq \bj(S)$ such that $\inv_\hfam(\pi)=S$. By the previous lemma $\pi|_{\bj(S)} \in \Hinv(S,\,\bj(S))$, which proves $\hinv(S;\bj(S))> 0$.
\end{proof}

By Lemma~\ref{lemma.flatten}, every permutation with $\hfam$-inversions equal to $S$ flattens to an element of $\Hinv(S, \bj(S))$. The next result studies the partition of $\Hinv(n,S)$ defined by the fibers of the flattening map. The key fact is that we can enumerate each such fiber using a polynomial formula. 

\begin{proposition}\label{prop.polynomial} Let $S$ be an $\hfam$-admissible set. For each $n\geq \bj(S)$, we have
\begin{eqnarray}\label{eqn.polyproofexpansion}
\hinv(S;n) = \sum_{\sigma\in \Hinv(S,\,\bj(S))} {n-\bt(\sigma)\choose \bj(S)-\bt(\sigma)}
\end{eqnarray}
where for each $\sigma\in \Hinv(S,\,\bj(S))$ we define  $\bt(\sigma) := \max\{\sigma_k \mid 
(k, \bj(S)+1)\in \cp_\hfam\}$.
\end{proposition}

\begin{proof} We prove the formula~\eqref{eqn.polyproofexpansion} using the flattening map. For each $\sigma\in \Hinv(S,\,\bj(S))$ consider the set
\[
\Hinv^{\sigma}(S,n):= \{\pi \in \Hinv(S,n) \mid \pi|_{\bj(S)} = \sigma\}.
\]
In other words, $\Hinv^\sigma(S,n)$ is the fiber of the flattening map $\Hinv(S,n) \to \Hinv(S,\,\bj(S))$ over $\sigma$. 
By Lemma~\ref{lemma.flatten}, these sets partition $\Hinv(S,n)$. 

Notice that, since $\sigma$ is a permutation of $[\bj(S)]$, we have $\bt(\sigma)\leq \bj(S)$.  To complete the proof, we claim  
\begin{eqnarray}\label{eqn.polyproofclaim}
\# \Hinv^{\sigma}(S,n) = {n-\bt(\sigma) \choose n-\bj(S)}.
\end{eqnarray}
Given this claim, as ${n-\bt(\sigma) \choose n-\bj(S)}= {n-\bt(\sigma)\choose \bj(S)-\bt(\sigma)}$, we recover the formula~\eqref{eqn.polyproofexpansion}.

For the rest of the proof, fix $\sigma\in \Hinv(S,\,\bj(S))$.
Let  $\pi \in \Hinv^{\sigma} (S,n)$. Consider the set $\{\pi_{\bj(S)+1}, \pi_{\bj(S)+2}, \ldots, \pi_n \}$. Recall that $\bm(S)$ denotes the maximum descent in $S$. Since $\bm(S)+1 \leq \bj(S)$, we know 
\[
\pi_{\bj(S)+1}<\pi_{\bj(S)+2}<\cdots < \pi_n.
\]
Consider a pair $(k, \bj(S)+1)$ in $\cp_\hfam$. Since $\bj(S)$ is the maximum index appearing in $S$ we have $(k, \bj(S)+1)\notin S$. Thus $\pi_k < \pi_{\bj(S)+1}$, as $\inv_\hfam(\pi)=S$. Furthermore, since $\pi|_{\bj(S)}=\sigma$ we have $\sigma_k \leq \pi_k$ for all $k\leq \bj(S)$. In particular, the fact that $\pi_k <\pi_{\bj(S)+1}$ for all $(k, \bj(S)+1)\in \cp_\hfam$ implies $\sigma_k < \pi_{\bj(S)+1}$ for all $(k, \bj(S)+1)\in \cp_\hfam$ as well. Because $\bt(\sigma)$ is the largest of these values $\sigma_k$, we conclude that $\bt(\sigma) < \pi_{\bj(S)+1}$. Thus, $\{\pi_{\bj(S)+1}, \pi_{\bj(S)+2}, \ldots, \pi_n\}$ is a cardinality $n-\bj(S)$ subset of  $[\bt(\sigma)+1, n] = \{\bt(\sigma)+1, \bt(\sigma)+2, \ldots, n\}$.

Let ${[\bt(\sigma)+1, n] \choose n-\bj(S)}$ denote the set of all cardinality $n-\bj(S)$ subsets of $[\bt(\sigma)+1, n]$. 
By the paragraph above, we have a well-defined map, 
\begin{eqnarray}\label{eqn.polyproofbijection}
\Hinv^\sigma(S,n) \to {[\bt(\sigma)+1, n] \choose n-\bj(S)},\;\; \pi \mapsto \{\pi_{\bj(S)+1}, \pi_{\bj(S)+2}, \ldots, \pi_n\}.
\end{eqnarray}
To prove our claim~\eqref{eqn.polyproofclaim}, we show this map is a bijection. The map is injective, since any permutation in $\Hinv(S,n)$ which flattens to $\sigma$ is uniquely determined by its elements in positions $\bj(S)+1, \bj(S)+2, \ldots, n$ (which are necessarily in increasing order).  

Next, we show that~\eqref{eqn.polyproofbijection} is surjective. Let $T \subseteq [\bt(\sigma)+1, n]$ be a subset of cardinality $n-\bj(S)$. Let $\pi$ be the permutation in $\Symm_n$ obtained by listing the elements of $[n]\setminus T$ in the same relative order as $\sigma$ and then listing the elements of $T$ in increasing order. By construction, $\pi|_{\bj(S)}= \sigma$. Since $1, 2, \ldots, \bt(\sigma)\in [n]\setminus T$, these values must appear in the same positions in $\pi$ as in $\sigma$. That is, 
\begin{eqnarray}\label{eqn.polyproof}
\pi_k = \sigma_k \; \text{ for all $k$ such that } \; \sigma_k \leq \bt(\sigma). 
\end{eqnarray}
Since $\pi$ flattens to $\sigma$ we have $S = \inv_\hfam (\sigma)\subseteq \inv_\hfam (\pi)$. Suppose, for a contradiction, that $(i,j) \in \inv_\hfam(\pi)\setminus S$. Then $(i,j)\in \cp_\hfam$ and $\pi_i>\pi_j$. Furthermore, since $(i,j)\notin S$, $\pi$ flattens to $\sigma$, and $\pi_{\bj(S)+1}<\pi_{\bj(S)+2}<\cdots < \pi_n$ we must have $i\leq \bj(S)$ and $\bj(S)+1\leq j$. The second inequality implies $\pi_j\in T$.  Since $(i,j)\in \cp_\hfam$ we get  $(i,\bj(S)+1)\in \cp_\hfam$ and so by definition of $\bt(\sigma)$, we conclude $\sigma_i \leq \bt(\sigma)$. It follows that $\pi_i = \sigma_i$ by~\eqref{eqn.polyproof}. Now the fact that $\pi_i > \pi_j$ implies $\sigma_i = \pi_i>\pi_j > \bt(\sigma)$ since $\pi_j\in T$, a contradiction. We conclude $\inv_\hfam(\pi)=S$, as desired. The map~\eqref{eqn.polyproofbijection} sends $\pi$ to $T$, so it is surjective. 
\end{proof}

Proposition~\ref{prop.polynomial} implies there is a polynomial expression in $n$ which computes the value of $\hinv(S;n)$ whenever it is nonzero.

Suppose $S$ is an $\hfam$-admissible set. In Proposition~\ref{prop.polynomial} above, we considered the maximum value $\bt(\sigma)$ obtained by a permutation $\sigma\in \Hinv(S,\bj(S))$ in positions $k \in [\bj(S)]$ such that $\hfam(k)>\bj(S)$. While these maximum values change as we consider different permutations $\sigma$ in the set $\Hinv(S,\bj(S))$, the next lemma tells us that the position where this maximum value occurs is uniquely determined.

\begin{lemma}\label{lem.t(S)existence} Let $S$ be an $\hfam$-admissible set, and $\bt(\sigma)$ be as defined in Proposition~\ref{prop.polynomial}. For all $\sigma,\pi\in \Hinv(S,\bj(S))$, the position in the one-line notation of $\sigma$ where $\bt(\sigma)$ occurs is equal to the position in the one-line notation of $\pi$ where $\bt(\pi)$ occurs.
\end{lemma}
\begin{proof} Let $p$ be the position in $\sigma$ where $\bt(\sigma)$ occurs and $q$ be the position in $\pi$ where $\bt(\pi)$ occurs, i.e., $p=\sigma^{-1}_{\bt(\sigma)}$ and $q = \pi^{-1}_{\bt(\pi)}$. By definition of the numbers $\bt(\sigma)$ and $\bt(\pi)$, $p,q\in [\bj(S)]$ such that $(p,\bj(S)+1), (q, \bj(S)+1)\in \cp_\hfam$, that is, such that $\hfam(p)>\bj(S)$ and $\hfam(q)>\bj(S)$. 

Assume $p\neq q$ and, without loss of generality, that $p<q$.
Since $\bt(\sigma) = \sigma_p > \sigma_q$, we know $(p,q)$ is an inversion of $\sigma$. Furthermore, we have $q\leq \bj(S)<\hfam(p)$ so $(p,q)\in \cp_\hfam$ and now the fact that $S=\inv_\hfam(\sigma)$ implies $(p,q)\in S$. 

On the other hand, $\bt(\pi) = \pi_q>\pi_p $ so $(p,q)$ is not an inversion of $\pi$. But this implies $(p,q) \notin \inv_\hfam(\pi)=S$, contradicting the conclusion above. Thus, $p$ and $q$ must be equal.
\end{proof}

Lemma \ref{lem.t(S)existence} tells us the number 
\[
\bt(S):= \sigma^{-1}_{\bt(\sigma)}\text{  for  }\sigma\in \Hinv(S,\bj(S))
\] 
is well-defined.

\begin{corollary} \label{cor.polynomial} Suppose $\hfam$ is a sequence of positive integers such that $\hfam(i)>i$ for all $i$ and let $S$ be an $\hfam$-admissible set. There exists a polynomial in $n$ whose value is equal to $\hinv(S;n)$ for all $n\geq \bj(S)$. In particular, for 
each $n\geq \bj(S)$, we have
\begin{eqnarray}\label{eqn.t_expansion}
\hinv(S;n) = \sum_{k=1}^{\bj(S)} t_k(S) {n-k\choose \bj(S)-k}
\end{eqnarray}
where for each $1\leq k \leq \bj(S)$ we define  $t_k(S) :=  \#T_k(S,\bj(S))$ with $T_k(S,\,\bj(S)):= \{\sigma \in \Hinv(S,\,\bj(S)) \mid \sigma_{\bt(S)}=k\}$.
\end{corollary}
\begin{proof} By Lemma~\ref{lem.t(S)existence}, $\bt(\sigma) = \sigma_{\bt(S)}$ for all $\sigma \in \Hinv(S,\bj(S))$. Thus, the expansion above is equivalent to that of Proposition~\ref{prop.polynomial}. 
\end{proof}

\begin{definition} For each $\hfam$-admissible set $S$, we call the polynomial appearing on the RHS of~\eqref{eqn.t_expansion} the \emph{$\hfam$-inversion polynomial of $S$}. By a slight abuse of notation, we denote this polynomial by $\hinv(S;n)$.
\end{definition}

Note that the formula~\eqref{eqn.t_expansion} does not count permutations in the set $\Hinv(S,n)$ for $1\leq n<\bj(S)$. Indeed, we know this number is $0$ for each such $n$, but the polynomial function on the RHS of~\eqref{eqn.t_expansion} does not, in general, have roots at all $1\leq n < \bj(S)$.

The next two examples compute the expansion for $\hinv(S;n)$ given in Corollary~\ref{cor.polynomial}.

\begin{example}\label{ex.polynomialproof} Let $\hfam(i)=i+3$ for all $i$ and  $S= \{(3,4), (3,5), (3,6), (4,6), (5,6)\}$. Then $\bj(S)=6$ and
\[
\Hinv(S, 6) = \{ 126453, 136452, 236451 \}.
\]
We have $\bj(S)+1=7$ and the values $i$ satisfying $i<7 \leq \hfam(i)=i+3$ are $i= 4, 5, 6$. Therefore $\bt(\sigma) = \max\{\sigma_4, \sigma_5, \sigma_6\}=5$ for all $\sigma\in \Hinv(S,6)$. Note that this maximal value occurs in position $\bt(S)=5$ for each permutation, confirming the results of Lemma~\ref{lem.t(S)existence}. The expansion~\eqref{eqn.t_expansion} is 
\[
\hinv(S;n) = 3 \, {n-5 \choose 6-5} = 3(n-5).
\]
\end{example}

\begin{example} Let $\hfam(i)=i+2$ for all $i$ and  $S= \{(1,3), (2,3), (2,4)\}$. Then $\bj(S)=4$ and
\[
\Hinv(S, 4) = \{ 2413, 3412 \}.
\]
We have $\bj(S) = 4$ and the values $i\leq \bj(S)=4$ satisfying $\hfam(i)>\bj(S) = 4$ are $i=3,4$. 
Therefore $\bt(\sigma) = \max\{\sigma_3, \sigma_4\}$. We get $\bt(2413) = 3$ and $\bt(3412)=2$ and $\bt(S) = 4$.  The expansion~\eqref{eqn.t_expansion}~is 
\[
\hinv(S;n) = {n-2\choose 4-2} + {n-3 \choose 4-3} = {n-2\choose 2}+ {n-3 \choose 1}.
\]
\end{example}

\begin{remark}\label{rem.descentexp} Consider the descent polynomial case, where $\hfam(i)=i+1$ for all $i\in \Z_{>0}$ and $S$ is a fixed descent set. In this case, $\bj(S)=m+1$ and the value of $\bt(\sigma)$ for $\sigma\in \Hinv(S,m+1)$ is just $\bt(\sigma)=\sigma_{m+1}$ so $\bt(S)=m+1$. Since $\sigma_{m}>\sigma_{m+1}$ we get $\sigma_{m+1}\in \{1, 2, \ldots, m\}$ and the expansion~\eqref{eqn.t_expansion} from Corollary~\ref{cor.polynomial} becomes
\[
\hinv(S;n) = \sum_{k=1}^m t_k(S) {n-k \choose m-k+1}
\]
where $t_k(S)= \#\{\sigma\in \Hinv(S,n) \mid \sigma_{m+1}=k\}$. This is precisely the expansion considered in the proof of \cite[Proposition 3.1]{Bencs}.   
\end{remark}

\begin{remark} \label{rem.non-polynomial} If we drop the assumption of $\hfam(i)>i$, and instead consider $\hfam(i)\geq i$, then the count $\hinv(S;n)$ may not be a polynomial. Indeed, if $\hfam(i)=i$ for all $i$ then $\varnothing$ is the only $\hfam$-admissible set and $\hinv(\varnothing;n)=n!$ is not a polynomial. Generally, $\hinv(S;n)$ is not a polynomial whenever $\{i \mid h(i) = i\}$ is not bounded.
\end{remark}

We can use the expansion from Corollary~\ref{cor.polynomial} to compute the degree of  $\hinv(S;n)$ and characterize when it is constant.  

\begin{corollary}\label{cor:constant}
Let $S\subseteq \cp_h$ be an $\hfam$-admissible set. The polynomial  $\hinv(S;n)$ has degree $\deg \hinv(S;n) = \bj(S) - \min \{ \sigma_{\bt(S)} \mid \sigma \in \Hinv(S,\bj(S)) \}$. In particular, it is a constant polynomial if and only if $\bt(\sigma) = \bj(S)$ for all $\sigma\in \Hinv(S,\bj(S))$. 
\end{corollary}
\begin{proof} The statement follows directly from Corollary~\ref{cor.polynomial} since $\deg {n-i \choose \bj(S)-i} = \bj(S)-i$ and since $\sigma_{\bt(S)} = \bt(\sigma)$. 
\end{proof}

\begin{example}\label{ex:constant2} 
Consider $\hfam = (2,4,4,5,6,7,8,\cdots)$ and $S = \{(2,3)\}$ from Example \ref{ex:constant_pre}, where we saw that $\hinv(S;n) = 2$ and $\Hinv(S,3) = \{132, 231\}$. The indices $i$ with $\hfam(i)>3$ are $2$ and $3$; we see $\bt(132) = \bt(231) = 3=\bj(S)$ in this case, so $\hinv(S;n)$ is indeed constant by Corollary~\ref{cor:constant}.  
\end{example}

%%%%%%%%%%%%%%%%%%%%%%%%%%

\section{Two expansion formulas for \texorpdfstring{$\hfam$}{hfam} weakly increasing}  \label{sec.expansion}

In this section, we focus on the case where $\hfam$ is a weakly increasing sequence of positive integers such that $\hfam(i)>i$ for all $i$, that is, $\hfam$ is a Hessenberg sequence. We give two expansions (which, in general, differ from that in Corollary~\ref{cor.polynomial}) for the $\hfam$-inversion polynomial and identify the coefficients explicitly. 

Throughout this section, $S$ denotes an $\hfam$-admissible set. As before, we set $\bm(S)=\max\{i\mid (i,i+1) \in S\}$. When $S$ is fixed we write $m=\bm(S)$. 

\subsection{Expansion Formula 1: The \texorpdfstring{$(b_k(S))$}{bk(S)}-expansion} \label{sec.expansion1} 
We generalize the expansion for the descent polynomial of~\cite[Proposition 3.1]{Bencs} to the setting of $\hfam$-inversion polynomials. In particular, our main result (Corollary~\ref{cor.k-expansion} below) computes the expansion of $\hinv(S;n)$ in the binomial basis 
\[
\left\{ {n-\hfam(m)\choose 0},{n-\hfam(m)+1\choose 1}, \ldots, {n-\hfam(m)+m\choose m}  \right\}
\]
and identifies the coefficients as nonnegative integers.

Given an $\hfam$-admissible set $S$, for each $n\geq \hfam(m)$ we define
\[
B_k(S,n):= \{\pi\in \Hinv(S,n) \mid  \pi_{\hfam(m)} = k \}.
\]
The sets $B_k(S,n)$ have a nice inductive description:

\begin{lemma}\label{lemma.k-decomp} Let $\hfam$ be a Hessenberg sequence, $S$ be an $\hfam$-admissible set, and $m= \bm(S)$. For all $n\geq \hfam(m)$ we have $B_k(S,n)\neq \varnothing$ only if $ \hfam(m)-m \leq k\leq \hfam(m)$. Furthermore, for each such $k$ there exists a bijection
\begin{eqnarray}\label{eqn.k-decomp}
B_k(S,n) \to  B_{k}(S,\hfam(m)) \times {[k+1, n] \choose n-\hfam(m)}.
\end{eqnarray}
\end{lemma}
\begin{proof} Suppose \( B_k(S, n) \neq \varnothing \), and let \( \pi \in B_k(S, n) \). Since \( \pi_{\hfam(m)} = k \) and 
\begin{eqnarray*}
\pi_{m+1}<\pi_{m+2}<\cdots< \pi_{\hfam(m)-1}< \pi_{\hfam(m)} < \pi_{\hfam(m)+1} \cdots < \pi_n,
\end{eqnarray*}
the number $\pi_{\hfam(m)}=k$, must be greater than at least $\hfam(m)-m$ many values (namely, $\pi_{m+1}, \ldots, \pi_{\hfam(m)-1}$) and must also be smaller than at least $n-\hfam(m)$ many values (namely, $\pi_{\hfam(m)+1}, \ldots, \pi_n$). Thus, $\hfam(m)-m \leq k\leq \hfam(m).$

We denote the function from~\eqref{eqn.k-decomp} by $f$. We define $f$ by
\[
f(\pi) =  (\pi|_{\hfam(m)} , \{\pi_{\hfam(m)+1}, \pi_{\hfam(m)+2}, \ldots, \pi_n\}).
\]
By Lemma~\ref{lemma.flatten}, $\pi|_{\hfam(m)}\in \Hinv(S,\hfam(m))$ and since the $n-\hfam(m)$ many values 
\[
\pi_{\hfam(m)+1}, \pi_{\hfam(m)+2}, \ldots, \pi_n
\]
are all larger than $\pi_{\hfam(m)}=k$, the map $f$ is well-defined.

We will show that $f$ is a bijection. Suppose $f(\pi) = f(\pi')$ for $\pi, \pi' \in B_k(S,n)$. Then, $\pi|_{[\hfam(m)]} = \pi'|_{[\hfam(m)]}$ and $\{\pi_{\hfam(m)+1}, \ldots, \pi_n\} = \{\pi'_{\hfam(m)+1}, \ldots, \pi'_n\}$. Since the entries at positions $\hfam(m)+1$ through $n$ must be increasing, $\pi$ and $\pi'$ must have the same entries in those positions. Similarly, since the first $\hfam(m)$ positions must appear in the same relative order and the remaining entries are fixed, there is a unique way to reconstruct the full permutation. Therefore, $\pi = \pi'$ and $f$ is injective.

Now, consider a permutation $\sigma\in B_k(S,\hfam(m))$ and a $n-\hfam(m)$ element subset $T$ of $[k+1, n]$. Define a permutation $\pi \in \Symm_n$ by listing the elements of $[n]\setminus T$ in the same relative order as $\sigma$ and then appending the elements of $T$ in increasing order. By construction, $\pi|_{\hfam(m)} = \sigma$. Since $1, 2, \ldots, k\in [n]\setminus T$, these values must appear in the same positions in $\pi$ as in $\sigma$. In particular, we know that $\sigma_{\hfam(m)}=k$ implies that $\pi_{\hfam(m)}=k$ as well. Furthermore, we know $\sigma_{m+1}<\sigma_{m+2}<\cdots <\sigma_{\hfam(m)-1}< \sigma_{\hfam(m)}=k$ so the values of $\pi_{m+1}, \pi_{m+2}, \ldots, \pi_{\hfam(m)}$ are all less than or equal to $k$. 

To complete the proof, we must show $\inv_\hfam(\pi)=S$.  Since $\pi$ flattens to $\sigma$ we have $S=\inv_\hfam(\sigma) \subseteq \inv_\hfam(\pi)$. To prove equality, suppose $(i,j)\in \inv(\pi)\setminus S$.  Then we have $i\leq \hfam(m)<j$ and $\pi_i>\pi_j$. Since $\pi_j\in T$ we conclude $\pi_i>k$. Since the values of $\pi_{m+1}, \ldots, \pi_{\hfam(m)}$ are all at most $k$ by the reasoning in the previous paragraph, we conclude $i\leq m$. Because $\hfam$ is increasing we know that $\hfam(i)\leq \hfam(m)$. Furthermore, $j>\hfam(m)$ by assumption. This implies that if $(i,j)\in \inv(\pi)\setminus (S)$ then $(i,j)\notin \cp_\hfam$ and therefore $(i,j)\notin \inv_\hfam(\pi)$. Thus  $\inv_\hfam(\pi)=S$, as desired. We conclude $\pi \in B_k(S,n)$ and since $f(\pi) = (\sigma, T)$, $f$ is indeed surjective.
\end{proof}

We obtain our expansion of the $\hfam$-inversion polynomial of $S$ as an immediate corollary.

\begin{corollary} \label{cor.k-expansion}  Let $\hfam$ be a Hessenberg sequence, $S$ be an $\hfam$-admissible set, and $m= \bm(S)$. Then
\begin{align*}
\hinv(S;n) &= \sum_{k=\hfam(m)-m}^{\hfam(m)} b_k(S) {n-k \choose \hfam(m)-k} \\
&= \sum_{\ell =0}^m b_{\hfam(m)-\ell}(S) {n-\hfam(m)+\ell \choose \ell}.
\end{align*}
where $b_k(S):=  \# B_k(S,\hfam(m))$.
\end{corollary}
\begin{proof} This is a direct consequence of Lemma~\ref{lemma.k-decomp}. Indeed,
the $B_k(S,n)$ sets form a partition of the set $\Hinv(S,n)$ and by Lemma~\ref{lemma.k-decomp}, $\# B_k(S,\hfam(m)) = 0$ unless $ \hfam(m)-m \leq k\leq \hfam(m)$. Furthermore, the bijection~\eqref{eqn.k-decomp} gives us
\begin{align*}
\# B_k(S,n) =  \# B_k(S,\hfam(m)) {n-k \choose n-\hfam(m)} =  \# B_k(S,\hfam(m)) {n-k \choose \hfam(m)-k}. 
\end{align*} 
The final equality follows by shifting the indexing of the sum, from $k$ to $\ell = \hfam(m)-k$.
\end{proof}

\begin{example}\label{ex.expb} Let $\hfam(i)=i+2$ for all $i\in \Z_{>0}$ and consider the $\hfam$-admissible set $S=\{(1,3), (2,3), (2,4)\}$ with $\bm(S)=2$. We have $\hfam(2)=4$ in this case and 
\[
\Hinv(S,\hfam(m)) = \Hinv(S,4) = \{2413, 3412\}.
\]
The coefficient $b_k(S)$ counts the number of permutations in $\Hinv(S,4)$ with entry in the last position equal to $k$. Therefore $b_2(S)=b_3(S)=1$. We conclude
\[
\hinv(S;n) = {n-2\choose 2} + {n-3\choose 1}.
\]
\end{example}

\begin{example} \label{ex.exp2} Let $\hfam(i)=i+3$ for all $i$ and  $S= \{(3,4), (3,5), (3,6), (4,6), (5,6)\}$. We have $\bm(S)=5$ and $\hfam(5)=8$ and in this case,
\begin{align*}
\Hinv(S,8) =  \{ 12645378, &12745368, 12845367, 13645278, 13745268, \\& \quad\quad\quad13845267, 23645178, 23745168, 23845167\}.
\end{align*}
The coefficient $b_k(S)$ counts the number of permutations in $\Hinv(S,8)$ with entry in the last position equal to $k$. Therefore $b_7(S)= 3$ and $b_8(S)=6$ in this case. We conclude
\[
\hinv(S;n) = 3 {n-7\choose 8-7} +  6 {n-8\choose 8-8} = 3(n-7) + 6 = 3n - 15,
\]
which agrees with Example~\ref{ex.polynomialproof}.
\end{example}

\subsection{Expansion Formula 2: The \texorpdfstring{$(a_k(S))$}{ak(S)}-expansion} We now prove a formula which generalizes the expansion formula \cite[Theorem 3.3]{descent_poly} of Diaz-Lopez, Harris, Insko, Omar, and
Sagan to the setting of restricted inversion polynomials.

\begin{theorem}\label{thm.a-expansion} Let $\hfam$ be a Hessenberg sequence, $S$ be an $\hfam$-admissible set, and $m= \bm(S)$. Then 
\begin{eqnarray}\label{eqn.a-expansion}
\hinv(S;n) = \sum_{k=0}^m a_k(S) {n-\hfam(m)+1\choose k}
\end{eqnarray}
where 
\[
a_k(S):= \#\{ \pi \in \Hinv(S,m+\hfam(m) -1) \mid \pi[m] = [\hfam(m), \hfam(m)+k-1]\}
\]
with $\pi[m]:= \{ \pi_1, \pi_2, \ldots, \pi_m\} \cap [\hfam(m), n]$.
\end{theorem}

\begin{proof} Define
\[
A_k (S,n) = \{\pi\in \Hinv(S,n) \mid \#\pi[m] = k \}.
\]
Since $\Hinv(S,n) = \bigsqcup_{k=0}^m A_k(S,n)$, it suffices to show
\begin{eqnarray}\label{eqn.a-expansion.claim}
\#A_k(S,n) = a_k(S) {n-\hfam(m)+1\choose k} \;\; \text{ for all } \;\; 0\leq k \leq m.
\end{eqnarray}
We may assume without loss of generality that $n\geq m+\hfam(m) -1$. Indeed, if we show that two polynomials on either side of~\eqref{eqn.a-expansion} are equal at infinitely many values, then they must be equal everywhere.

First, consider the set 
\[
A_k^*(S,n) = \{\pi \in A_k(S,n) \mid \pi[m] = [\hfam(m), \hfam(m)+k-1]\}.
\]
By definition, $a_k(S) = \# A_{k}^*(S, m + \hfam(m)-1)$. 
Consider the map 
\[
A_k^*(S, m + \hfam(m) - 1) \to A_k^*(S,n)
\]
defined by taking a permutation \( \pi \) on the left-hand side and appending the elements in \( [m+ \hfam(m), n] \) in increasing order to obtain a permutation in \( A_k^*(S, n) \).

Let \( \pi \in A_k^*(S, n) \). Since \( (m, m+1) \) is the last descent of $\pi$, we must have 
\[
\pi_{m+1}<\pi_{m+2}< \cdots <\pi_{m+\hfam(m)-1} < \pi_{m+\hfam(m)}<\cdots <\pi_n.
\]  
Thus, $\pi_{i}\leq m+\hfam(m)-1$ for all $i\in [m+1, m+\hfam(m)-1]$. Furthermore, since \( \pi[m] = [\hfam(m), \hfam(m)+k-1] \) and $k\leq m$, we have 
\[
\pi_i<\hfam(m)+k-1 \leq \hfam(m)+m-1
\]
for all $i\in [m]$.  We conclude  \( \pi_1, \ldots, \pi_{m+\hfam(m)-1} \in [1, m+\hfam(m)-1] \). This forces $\pi_{i}=i$ for all $i\in [m+\hfam(m), n]$, proving our map from $A_k^*(S,m+\hfam(m)-1)$ to $A_k^*(S,n)$ is surjective. Since the map is also clearly injective, we have
\begin{eqnarray}\label{eqn.a-expansion0}
\# A_k^*(S,n) = \# A_k^*(S,m+\hfam(m)-1) = a_k(S)
\end{eqnarray}
for all $0\leq k \leq m$. 

Note that when $k=0$, we get $A_0^*(S,n) = A_0(S,n)$ so~\eqref{eqn.a-expansion.claim} holds in this case.

Now suppose that $k\geq 1$.  Let $X$ and $Y$ be subsets of $[\hfam(m), n]$ of cardinality $k$ such that the sets
\[
 \{\pi\in A_k(S,n) \mid \pi[m]=X\} \; \text{ and } \;  \{\pi\in A_k(S,n) \mid \pi[m]=Y\}
\]
are nonempty. 
There is a unique, order preserving bijection $f: X \to Y$. We use $f$ to define a map 
\begin{equation}\label{eq.a_expansionk}
    \tilde{f}: \{\pi\in A_k(S,n) \mid \pi[m]=X\} \to \{\pi\in A_k(S,n) \mid \pi[m]=Y\}
\end{equation}
as follows. Suppose $\pi\in A_k(S,n)$ such that $\pi[m]=X$. We obtain the permutation $\tilde{f}(\pi)$ from $\pi$ by applying $f$ to the elements of $\pi[m]$ in the one-line notation of $\pi$, leaving all entries from $[\hfam(m)-1]$ unchanged. We then list all remaining entries of $[\hfam(m), n]\setminus Y$ in increasing order to obtain the permutation $\tilde{f}(\pi)$. 

To prove $\tilde{f}$ is a bijection, it suffices to show $\tilde{f}(\pi)\in A_k(S,n)$. Indeed, given this fact we obtain the inverse of $\tilde{f}$ by applying the same construction to $f^{-1}: Y \to X$. Let $\sigma = \tilde{f}(\pi)$. We aim to show that $\inv_\hfam(\sigma) = S$, or in other words, that $\inv_\hfam(\pi) = \inv_\hfam(\sigma)$.

First suppose $(i,j)\in \inv_\hfam(\pi)$. We consider three possible cases:

\

\textbf{Case 1:} Assume $i,j\leq m$. Since $(i,j)$ is an inversion of $\pi$ we must have $\pi_i>\pi_j$. There are three possible subcases:
\begin{itemize}
\item $\pi_i, \pi_j<\hfam(m)$, so $\sigma_i=\pi_i$ and $\sigma_j = \pi_j$. Then $\sigma_i = \pi_i > \pi_j = \sigma_j$ so $(i,j)\in \inv_{\hfam}(\sigma)$. 
\item $\pi_j<\hfam(m)$ and $\pi_i \geq \hfam(m)$, so $\sigma_i = f(\pi_i)$ and $\sigma_j = \pi_j$. Then $\sigma_i = f(\pi_i) \geq \hfam(m) > \pi_j = \sigma_j$ so $(i,j)\in \inv_\hfam(\sigma)$. 
\item $\pi_i, \pi_j>\hfam(m)$, so $\sigma_i= f(\pi_i)$ and $\sigma_j = f(\pi_j)$. Then $(i,j)\in \inv_\hfam(\sigma)$ since $f$ is an order preserving bijection. 
\end{itemize}
In each of these subcases, we have confirmed that $(i,j)\in \inv_\hfam(\sigma)$, as desired. 

\

\textbf{Case 2:} Assume $i,j>m$.  Since all the entries of $\pi$ in positions $m+1, m+2 , \ldots, n$ are in increasing order, $(i,j)$  is not an inversion of $\pi$ so this case never occurs.

\

\textbf{Case 3:} Assume $i\leq m$ and $j>m$. Since $(i,j)\in S$ we know $j\leq \hfam(i)$. Because $i\leq m$ this tells us that $j\leq \hfam(m)$ so $j\in [m+1, \hfam(m)]$. (Note that $\hfam(m)\geq m+1$ since $(m,m+1)\in S$, so this interval is well-defined.)  We want to show $\sigma_i>\sigma_j$. 

We know all entries of $\pi$ in positions $m+1, m+1 , \ldots, n$ are in increasing order. Since $\#X=k$, there are $m-k$ many numbers from the set $[\hfam(m)-1]$ that appear in the first $m$ entries of the one-line notation for $\pi$. The numbers from the set $[\hfam(m)-1]$ that do not appear in the first $m$ entries must be listed in order, starting in position $m+1$. There are 
\[
\#\left( [\hfam(m)-1] \setminus \{\pi_1,\pi_2, \ldots, \pi_m\}  \right) = \hfam(m)-1 - (m-k) \geq \hfam(m)-m
\]
many such values, since $k\geq 1$.  In particular, it follows that the entries of $\pi$ in positions $m+1, m+2, \ldots, \hfam(m)$ are elements of $[\hfam(m)-1]$. Since $j\in [m+1, \hfam(m)]$, we have $\pi_j \in [\hfam(m)-1]$ and thus $\sigma_j=\pi_j$. 

Now, if $i\in [\hfam(m)-1]$ also, then $\sigma_i = \pi_i$ so $\sigma_i=\pi_i>\pi_j = \sigma_j$ as desired. Otherwise, $i\geq \hfam(m)$ and $\sigma_i = f(\pi_i) \geq \hfam(m) > \hfam(m)-1 \geq \pi_j=\sigma_j$. Thus $(i,j)\in \inv_\hfam(\sigma)$.

\

We have now proved that $\inv_\hfam(\pi)\subseteq \inv_\hfam(\sigma)$. The proof of the opposite inclusion follows in exactly the same way, so we leave it to the reader.

We now conclude 
\begin{align*}
\# A_k(S,n) &= \bigsqcup_{\substack{X\subseteq [\hfam(m), n]\\ \#X=k}} \#\{\pi \in A_k(S,n) \mid \pi[m]=X\} \\
&= \# A_k^*(S,n) \cdot {n-\hfam(m)+1\choose k} \quad \text{ (because of the bijection~\eqref{eq.a_expansionk})} \\
&= a_k(S) {n-\hfam(m)+1\choose k} \quad \text{(by~\eqref{eqn.a-expansion0})},
\end{align*}
which completes the proof of~\eqref{eqn.a-expansion.claim} for $1\leq k \leq m$. 
\end{proof}

\begin{example}\label{ex.a-expansion} Let $\hfam(i)=i+2$ for all $i\in \Z_{>0}$ and $S=\{(1,3), (2,3), (2,4)\}$. In this case $m=\bm(S)=2$ and $m+\hfam(m)-1 = 2+4-1=5$.
By definition
\begin{align*}
a_0(S) &=\#\{\pi \in \Hinv(S, 5) \mid \{\pi_1, \pi_2\}\cap \{4,5\} = \varnothing \}\\
a_1(S) &=\#\{\pi \in \Hinv(S, 5) \mid \{\pi_1, \pi_2\}\cap \{4,5\} = \{4\} \}\\
a_2(S) &=\#\{\pi \in \Hinv(S, 5) \mid \{\pi_1, \pi_2\}\cap \{4,5\} = \{4,5\} \}.
\end{align*}
Since 
\[
\Hinv(S, 5) = \{24135, 25134, 34125, 35124, 45123\}
\]
we get $a_0(S)=0$, $a_1(S) = 2$ and $a_2(S)=1$. The expansion from Theorem~\ref{thm.a-expansion} gives us,
\[
\hinv (S;n) = 2 {n-3\choose 1} + {n-3 \choose 2} = 2(n-3) + \frac{1}{2}(n-3)(n-4).
\]
The reader can confirm that, although this expansion differs from that in Example~\ref{ex.expb}, the two formulas are indeed equal. 
\end{example}

\begin{lemma}\label{lem:ab_conversion} Suppose $S$ is a nonempty $\hfam$-admissible set. 
The coefficients $(a_k(S))_{k=0}^m$ and $(b_k(S))_{k=0}^m$ are related as follows:
\begin{eqnarray}\label{eqn.coeff}
a_0(S) = b_{\hfam(m)}(S) \; \text{ and } \; a_k(S) = \sum_{j=k}^m {j-1 \choose k-1} b_{\hfam(m)-j} \;\; \text{ for all $1\leq k \leq m$.}
\end{eqnarray}
\end{lemma}
\begin{proof} Using Vandermonde's identity for binomial coefficients we get that 
\begin{align*}
\hinv(S;n) &= \sum_{\ell=0}^m b_{\hfam(m)-\ell}(S) {n-\hfam(m)+\ell \choose \ell}\\
&= b_{\hfam(m)-0}(S) {n-\hfam(m)+0 \choose 0} + \sum_{\ell=1}^m b_{\hfam(m)-\ell}(S) {n-\hfam(m)+1 + (\ell-1) \choose \ell}\\
&= b_{\hfam(m)}(S) + \sum_{\ell=1}^m b_{\hfam(m)-\ell}(S) \left( \sum_{k=0}^\ell {n-\hfam(m)+1\choose k}{\ell-1\choose \ell-k} \right)
\end{align*}
Now if $k = 0$ or $k > \ell$ then ${\ell - 1 \choose \ell - k} = 0$, so we may let the second sum be from $k=1$ to $m$ instead of $k=0$ to $\ell$. 
Re-writing the sum above gives us
\[
\hinv(S;n) = b_{\hfam(m)}(S) + \sum_{k=1}^m \left( \sum_{\ell=1}^m {\ell-1 \choose \ell-k} b_{\hfam(m)-\ell}(S)  \right)   {n-\hfam(m)+1\choose k}.
\]
Comparing with the sum formula of Theorem~\ref{thm.a-expansion} and noting ${\ell-1 \choose \ell-k} = {\ell-1 \choose k-1} = 0$ whenever $\ell < k$ gives us~\eqref{eqn.coeff}.
\end{proof}

% \begin{lemma}\label{lem:at_conversion}
% Suppose $S$ is a nonempty $\hfam$-admissible set. 
% The coefficients $(a_k(S))_{k=0}^m$ and $(t_k(S))_{k=0}^{\bj(S)}$ are related as follows:
% \begin{eqnarray}\label{eqn.atcoeff}
% a_k(S) = \sum_{i=0}^{\bj(S)} {\hfam(m)-i-1 \choose \bj(S)-i-k}t_i(S)\;\; \text{ for all $0\leq k \leq m$.}
% \end{eqnarray}
% \end{lemma}
% \begin{proof}
%     Using Vandermonde's identity for binomial coefficients we get that 
%     \begin{align*}
%         \hinv(S,n) &= \sum_{i=0}^{\bj(S)} t_i(S) {n-i \choose \bj(S)-i} \\
%           &= \sum_{i=0}^{\bj(S)} t_i(S) {n-\hfam(m)+1+\hfam(m)-i-1 \choose \bj(S)-i} \\
%           &= \sum_{i=0}^{\bj(S)}  t_i(S)\sum_{k=0}^{\bj(S)-i}{\hfam(m)-i-1 \choose \bj(S)-i-k}{n-\hfam(m)+1 \choose k} \\
%     \end{align*}
%     Now if $k > \bj(S) - i$ then ${m-i-1 \choose \bj(S)-i-k} = 0$, so we may take the second sum be over all $k \geq 0$, allowing us to switch the sums and write 
%     \[
%     \hinv(S;n) = \sum_{k\geq 0}\left(\sum_{i=0}^{\bj(S)} t_i(S){\hfam(m)-i-1 \choose \bj(S)-i-k}\right){n-\hfam(m)+1 \choose k}.
%     \]
% Comparing with the expansion formula of Theorem~\ref{thm.a-expansion} gives us \eqref{eqn.atcoeff}.
% \end{proof}

% \begin{remark}
%     We note that Lemma \ref{lem:at_conversion} is perfectly valid when $\hfam$ is not weakly increasing, however we have no guarantee that the highest term in the $(a_k(S))$ expansion is $k=m$, nor that each $\hfam(m)-i-1$ is non-negative.
% \end{remark}

When $\hfam$ is not a weakly increasing sequence, the coefficient sequences $(a_k(S))$ and $(b_k(S))$ of $\hinv(S;n)$ in the relevant bases need not be positive. 

\begin{example}\label{ex:ab_negative}
Consider $\hfam = (7,3,4,5,7,8,9,10,\ldots)$ and 
\[
S = \{(1,4),(1,5),(1,6),(2,3),(3,4)\}.
\]
We get that $\bj(S) = 6$ and $\bt(S) = 1$ in this case. Applying Corollary \ref{cor.polynomial} using 
\[
\Hinv(S,6) = \{465123\},
\]
we find $t_1(S)=t_2(S)=t_3(S)=t_5(S)=t_6(S)=0$ and $t_4(S)=1$, thus 
\[
\hinv(S;n) = {n-4 \choose 2}.
\]
Note that $\hfam(m)-1 = \hfam(3)-1= 4-1=3$. Using Vandermonde's identity gives
\begin{align*}
    \sum_{k \geq 0}a_k(S) {n-\hfam(m)+1 \choose k}
     %&= {4-4-1 \choose 6-4-0}{n-3 \choose 0} + {4-4-1 \choose 6-4-1}{n-3 \choose 1} + {4-4-1 \choose 6-4-2}{n-3 \choose 2} \\
     &= {-1 \choose 2}{n-3 \choose 0} + {-1 \choose 1}{n-3 \choose 1} + {-1 \choose 0}{n-3 \choose 2}\\
     &= {n-3 \choose 0} -{n-3 \choose 1} + {n-3 \choose 2}
\end{align*} 
and hand calculations confirm 
\[\sum_{\ell \geq 0}b_{n-\ell}(S){n-\hfam(m)+\ell \choose \ell} ={n-4 \choose 0} - 2{n-3 \choose 1} +  {n-2 \choose 2}.\]
Thus, we are not guaranteed non-negative coefficients for the $(a_k(S))$ or $(b_k(S))$ expansions of $\hinv(S;n)$ when $\hfam$ is not weakly increasing.
\end{example}

\section{Log-concavity}
This section considers the nonnegative coefficient sequence $(t_k(S))_{k=1}^{\bj(S)}$ associated to an $\hfam$-admissible set $S$ and also the sequences $(a_k(S))_{k=0}^m$ and $(b_k(S))_{k=\hfam(m)-m}^{\hfam(m)}$ for an $\hfam$-admissible set $S$ when $\hfam$ is a Hessenberg sequence. We prove these sequences are log-concave, and in doing so, we will also prove that these sequences do not have internal zeros. The results in this section directly generalize the work of Bencs in~\cite[\S 3]{Bencs}.

Recall that we say the finite sequence $a = (a_0,a_1,a_2,\ldots,a_{\ell})$ of real numbers is \emph{log concave} if $(a_i)^2 \geq a_{i-1}a_{i+1}$ for all $i\in  [\ell-1]$. We say that $a$ is \emph{unimodal} if it first weakly increases and then weakly decreases, or more formally, if there is some $i$ such that 
\[
a_0 \leq \cdots \leq a_{i-1} \leq a_i \geq a_{i+1} \cdots \geq a_\ell.
\] 

\begin{remark} Suppose the sequence $a = (a_0,a_1,a_2,\ldots,a_{\ell})$ is nonnegative and that for all $a_i,a_k \neq 0$, if $i<j<k$ then $a_j \neq 0$ (i.e., $a$ has no \emph{internal zeros}). Then $a$ being log concave implies that it is also unimodal. This is a standard approach to establishing unimodality, and the unimodality of various integer sequences is an active area of research in algebraic combinatorics.
\end{remark}

We show that the coefficient sequences $(t_k(S))_{k=1}^{\bj(S)}$ and $(b_k(S))_{k=\hfam(m)-m}^{\hfam(m)}$ are log-concave by providing a different combinatorial interpretation for them, and applying a theorem of Stanley \cite{STANLEY1981}. We begin with some background on the height polynomials of partially ordered sets; see~\cite{ChungFishburnGraham_heightpoly,STANLEY1981,Bencs}.

\subsection{Height polynomials}
Let $P$ be a finite poset on $n$ elements and let $v \in P$ be a fixed element. The \textit{height polynomial} $h_{P,v}(x)$ is defined by:
\[
h_{P,v}(x) = \sum_{\phi \in \text{Ext}(P)} x^{\phi(v)-1} = \sum_{k=0}^{n-1} h_k(P,v) x^k,
\]
where $\mathrm{Ext}(P)$ denotes the set of all linear extensions of $P$, and $\phi(v)$ is the position of $v$ in the linear extension $\phi\in \mathrm{Ext}(P)$. So, the number $h_k(P,v)$ counts the number of linear extensions in which exactly $k$ elements are less than $v$. We call $\left(h_k(P,v) \right)_{k=0}^{n-1}$ the \emph{height sequence of~$v$}. 

\begin{example}\label{ex:hpoly} Let $P$ be the poset on $[6]$ whose cover relations are $\{(1<_P 2),(2<_P3),(2<_P4),(3<_P5),(5<_P6),(4<_P6)\}$. The Hasse diagram of this poset is the following.
    
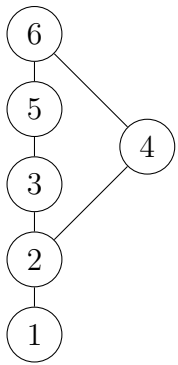
\begin{center}
    \begin{tikzpicture}
        \draw (0,0) node[draw,circle] (v1) {$1$};
        \draw (0,1) node[draw,circle] (v2) {$2$};
        \draw (0,2) node[draw,circle] (v3) {$3$};
        \draw (1.5,2.5) node[draw,circle] (v4) {$4$};
        \draw (0,3) node[draw,circle] (v5) {$5$};
        \draw (0,4) node[draw,circle] (v6) {$6$};

        \draw (v1)--(v2)--(v3)--(v5)--(v6);
        \draw (v2)--(v4)--(v6);
    \end{tikzpicture}
    \end{center}
The set of all linear extensions of this poset is precisely
    \[\mathrm{Ext}(P) = \{124356,\;\;123456,\;\;123546\}.\]
It follows that 
    \begin{center}
    \begin{tabular}{rclcrcl}
        $h_{P,1}(x)$ & $=$ & $3$ &{\hspace{1cm}} & $h_{P,4}(x)$ & $=$ & $x^2+x^3+x^4$\\
        $h_{P,2}(x)$ & $=$ & $3x$ &{} & $h_{P,5}(x)$ & $=$ & $x^3+2x^4$\\
        $h_{P,3}(x)$ & $=$ & $2x^2+x^3$ &{} & $h_{P,6}(x)$ & $=$ & $3x^5$.
    \end{tabular}
    \end{center}
    
\end{example}

The key result for our purposes is that Stanley proved that the coefficients of a height polynomial are log-concave~\cite[Thm 3.1]{STANLEY1981}.

\begin{theorem}[Stanley]\label{thm:heightpoly_logconcave}
Let $P$ be a finite poset on $n$ elements and $v \in P$. Then the height sequence $\left(h_k(P,v)\right)_{k=0}^{n-1}$ is log-concave. 
\end{theorem}
Theorem \ref{thm:heightpoly_logconcave} above is an application of a stronger result of Stanley, proved using the Alexandrov--Fenchel inequalities from the theory of mixed volumes. Lemma \ref{lem:height_internal_zeros} below also follows from Stanley's theorem, but also has a much simpler combinatorial proof, so we include it for completeness.

\begin{lemma}\label{lem:height_internal_zeros}
Let $P$ be a finite poset on $n$ elements and $v \in P$. The height sequence $\left(h_k(P,v)\right)_{k=0}^{n-1}$ has no internal zeros.
\end{lemma}

\begin{proof}
We write $h_k \coloneqq h_k(P,v)$ to simplify the notation in this proof. Let 
\[
P_\downarrow(v) \coloneqq \{w \in P \mid w <_P v\}\;\; \text{ and }\;\; P_\uparrow(v)\coloneqq \{w \in P \mid w >_P v\}
\]
be the set of elements in $P$ (strictly) less than and greater than $v$, respectively. \par

The first nonzero entry of the sequence $\left(h_k \right)_{k=0}^{n-1}$ must occur for index value at least $k = \#P_{\downarrow}(v)$. Indeed, in any linear extension of $P$, every element of $P_\downarrow(v)$ must appear before $v$. On the other hand, the last nonzero entry of the sequence $\left(h_k \right)_{k=0}^{n-1}$ must occur for index value at most $k = \# (P \setminus P_\uparrow(v))-1$. This is because every element of $P_\uparrow(v)$ must appear after $v$ in any linear extension of $P$. Thus $h_k=0$ for all $k < \#P_{\downarrow}(v)$ or $k > \# (P \setminus P_\uparrow(v))-1$. \par

Now given any $k$ such that $\#P_{\downarrow}(v) \leq k \leq \# (P \setminus P_\uparrow(v))-1$, we will construct a linear extension $\phi\in \mathrm{Ext}(P)$ such that $\phi(v)=k+1$.  
We let $\phi= \phi_1\phi_2\phi_3$ where $\phi_1$, $\phi_2$, and $\phi_3$ are constructed as follows.
\begin{enumerate}
\item Let $\phi_1$ be a linear extension of $P_\downarrow(v)$.
\item Let $\phi_2'$ be any linear extension of the set $P \setminus (P_\downarrow(v) \cup P_\uparrow(v) \cup \{v\})$. Note that this is the set of all elements of $P$ which are incomparable to $v$. The length of $\phi_2'$ is $n-\#P_\downarrow(v)-\#P_\uparrow(v)-1$. Insert $v$ in the $(k-\#P_\downarrow(v)+1)$-st position of $\phi_2'$ to construct $\phi_2$.
\item Let $\phi_3$ be a linear extension of $P_\uparrow(v)$.
\end{enumerate}
It's straightforward to check that $\phi$ satisfies the order relations for $P$. The linear extension $\phi\in \mathrm{Ext}(P)$ contributes a term of $x^k$ to the height polynomial. In other words, this shows that $h_k \geq 1$ for all $\#P_{\downarrow}(v) \leq k \leq \# (P \setminus P_\uparrow(v))-1$, so $\left(h_k\right)_{k=0}^{n-1}$ contains no internal zeros.
\end{proof}
\begin{example}\label{ex:hpoly_internal_zeros}
This example demonstrates the proof technique of Lemma~\ref{lem:height_internal_zeros}. Consider the poset $P$ from Example \ref{ex:hpoly} and let $v=4$. We have $\left(h_k(P,4)\right)_{k=0}^5 = (0,0,1,1,1,0)$ and
\[
P_\downarrow(4)=\{1,2\}, \;\; P_{\uparrow}(4)=\{6\},  \;\; \text{ and }\;\; [6]\setminus (P_\downarrow(4) \cup P_{\uparrow}(4) \cup\{4\}) = \{3,5\}.
\]
The last set is the set of all elements which are incomparable to $4$. The only linear extensions of these sets are $\phi_1 = 12$, $\phi_3 = 6$, and $\phi_2' = 35$, respectively. We can insert $4$ in any of the three places available in $\phi_2'$ to obtain $\phi_2$, and for each of these three possibilities, $\phi=\phi_1\phi_2\phi_3$ is one of the three linear extensions in $\mathrm{Ext}(P)$. Last, each of these possibilities corresponds to $4$ appearing in positions $3, 4, 5$, respectively, in $\phi$. This explains why $h_2(P,4) = h_3(P,4)=h_4(P,4)=1$. 
\end{example}

\subsection{Log-concavity for \texorpdfstring{$\left(t_k(S)\right)$}{tk(S)}}

We now connect our work on restricted inversion polynomials to height polynomials. 
Let $\hfam$ be a sequence of positive integers such that $\hfam(i)>i$ (not necessarily weakly increasing) and $S$ be an $\hfam$-admissible set.
For a positive integer $\ell\geq \bj(S)$ we define a poset $P_{\hfam, S}^{\ell}=([\ell], <_{S})$ on the set $[\ell]$ to be the transitive closure of the relation $<_S$ such that:
\[
i >_S j \; \text{ if }\; (i,j) \in S \quad \text{ and } \quad
i <_S j \; \text{ if }\; (i,j) \in S^c = \cp_\hfam\setminus S \; \text{ and }\; j\leq \ell.
\]
In this section we consider the case where $\ell = \bj(S)$.

\begin{example}\label{ex.partial_order} Consider the sequence 
\[
\hfam= (3,5,4,6,8,7,9,11,10,12,\ldots)
\]
and $\hfam$-admissible set $S=\{(2,4), (2,5), (3,4)\}$ with $\bj(S)=5$. 
In this case, the relations of $P_{\hfam, S}^{\bj(S)} = P_{\hfam, S}^{5}$ are the transitive closure of:
\begin{align*}
&1 <_S 2, \; 1<_S 3 ,\; 2<_S 3,\; 4<_S 5  &\qquad(\text{since } (1,2), (1,3), (2,3), (4,5) \in S^c)\\
&2 >_S 4,\; 2>_S 5, \;  3 >_S 4    &\qquad (\text{since } (2,4), (2,5), (3,4)\in S ).
\end{align*}
The Hasse diagram for $P_{\hfam, S}^{5}$ is as follows.
\begin{center}
    \begin{tikzpicture}
        \draw (1,0) node[draw,circle] (v1) {$4$};
        \draw (1,1) node[draw,circle] (v2) {$5$};
        \draw (0,1) node[draw,circle] (v3) {$1$};
        \draw (1,2) node[draw,circle] (v4) {$2$};
        \draw (1,3) node[draw,circle] (v5) {$3$};

        \draw (v1)--(v2);
        \draw (v3)--(v4);
        \draw (v2)--(v4);
        \draw (v4)--(v5);
    \end{tikzpicture}
\end{center}
Notice that the set of linear extensions of the poset is $\mathrm{Ext}(P_{\hfam, S}^{5}) = \{14523, 41523, 45123\}$. Comparing with the set
\[
\Hinv(S,\bj(S)) = \Hinv(S,5) = \{14523, 24513, 34512\},
\]
we see that the inverse of every permutation in $\Hinv(S,5)$ is a linear extension of $P_{\hfam, S}^5$.
\end{example}

Any linear extension of $P_{\hfam, S}^{\,\bj(S)}$ corresponds to a permutation in the set $\Hinv(S, \bj(S))$. In particular:

\begin{lemma}\label{lem:Tset_poset_map} Let $S$ be an $\hfam$-admissible set. Then 
\begin{eqnarray}\label{eqn.t-to-ht} 
t_k (S) = h_{k-1}(P_{\hfam, S}^{\,\bj(S)}, \bt(S)) \quad \text{ for all } \quad  1 \leq k \leq \bj(S). 
\end{eqnarray}
\end{lemma}
\begin{proof} By definition, $\pi \in \Hinv(S,\bj(S))$ if and only if $\pi^{-1}$ is a bijection from $[\bj(S)]$ to $[\bj(S)]$ that preserves the partial order $P_{\hfam, S}^{\,\bj(S)}$. Moreover, if $\pi \in \Hinv(S,\bj(S))$ such that $\pi_{\bt(S)}=k$, then there are precisely $k-1$ elements that precede $\bt(S)$ in $\pi^{-1}$.  Thus $h_{k-1}(P_{\hfam, S}^{\,\bj(S)},\bt(S))$ counts the permutations $\pi^{-1}$ such that $\pi_{\bt(S)} = k$, which is precisely $t_k(S)$.
\end{proof}

As an immediate consequence of Stanley's Theorem~\ref{thm:heightpoly_logconcave} and Lemmas~\ref{lem:height_internal_zeros} and~\ref{lem:Tset_poset_map}, we obtain the main result:

\begin{proposition}\label{prop:t_logconcave} Let $S$ be an $\hfam$-admissible set. The sequence of non-negative integers   $\left(t_k(S)\right)_{k=1}^{\bj(S)}$  is log-concave and does not have internal zeros. 
\end{proposition}

Using Lemma~\ref{lem:Tset_poset_map}, we can also rephrase Corollary \ref{cor:constant} in terms of the poset $P_{\hfam, S}^{\,\bj(S)}$:

\begin{corollary} Let $S$ be an $\hfam$-admissible set. The polynomial $\hinv(S;n)$ is constant if and only if $\bt(S)$ is the unique maximal element of $P_{\hfam, S}^{\,\bj(S)}$.
\end{corollary}

\begin{example} Consider $\hfam = (2,4,4,5,6,7,8,\cdots)$ and $S = \{(2,3)\}$ from Example~\ref{ex:constant_pre} and also Example~\ref{ex:constant2}, where we saw that $\hinv(S;n) = 2$ is constant and $\bt(S)=2$. The poset $P_{\hfam, S}^{\, \bj(S)} = P_{\hfam, S}^{3}$ has the following Hasse diagram, with unique maximal element equal to $\bt(S)=2$. 
\begin{center}
    \begin{tikzpicture}
        \draw (1.5,1) node[draw,circle] (v3) {$3$};
        \draw (0.5,1) node[draw,circle] (v2) {$1$};
        \draw (1,2) node[draw,circle] (v1) {$2$};

        \draw (v1)--(v2);
        \draw (v1)--(v3);
    \end{tikzpicture}
\end{center}
\end{example}

\subsection{Log-concavity for \texorpdfstring{$\left(b_k(S)\right)$}{bk(S)}} Now assume that $\hfam$ is a Hessenberg sequence and $S$ an $\hfam$-admissible set. The poset defined in the previous section gives a new interpretation of the $(b_k(S))$ coefficient sequence when $\ell = \hfam(m)$.
Any linear extension of $P_{\hfam, S}^{\hfam(m)}$ corresponds to a permutation in the set $\Hinv(S, \hfam(m))$. In particular:

\begin{lemma}\label{lem:Bset_poset_map} Let $S$ be an $\hfam$-admissible set and $m=\bm(S)$ denote the maximum descent in $S$. Then 
\begin{eqnarray}\label{eqn.b-to-ht} 
b_k (S) = h_{k-1}(P_{\hfam, S}^{\hfam(m)}, \hfam(m)) \quad \text{ for all } \quad  \hfam(m)-m \leq k \leq \hfam(m). 
\end{eqnarray}
\end{lemma}
\begin{proof} By definition, $\pi \in \Hinv(S,\hfam(m))$ if and only if $\pi^{-1}$ is a bijection from $[\hfam(m)]$ to $[\hfam(m)]$ that preserves the partial order $P_{\hfam, S}^{\hfam(m)}$. Moreover, if $\pi \in \Hinv(S,\hfam(m))$ such that $\pi_{\hfam(m)}=k$, then there are precisely $k-1$ elements that precede $\hfam(m)$ in $\pi^{-1}$.  Thus $h_{k-1}(P_{\hfam, S}^{\hfam(m)},\hfam(m))$ counts the permutations $\pi^{-1}$ such that $\pi\in B_k(S,\hfam(m))$.
\end{proof}

\begin{example}\label{ex.poset} This example demonstrates the equality of Lemma~\ref{lem:Bset_poset_map}.
Let $\hfam(i)=i+2$ for all $i$ and consider the $\hfam$-admissible set 
\[
S = \{(1,3),(2,3),(2,4),(3,4)\}. 
\]
with $\bm(S) =m = 3$.  In this case, the relations of $P_{\hfam, S}^{\hfam(m)} = P_{\hfam, S}^{5}$ are the transitive closure of:
\begin{align*}
&1 <_S 2, \; 3<_S5 ,\; 4<_S 5  \qquad & (\text{since } (1,2), (3,5), (4,5) \in S^c)\\
&1 >_S 3,\; 2>_S 3, \;  2 >_S 4, 3>_S4   \qquad &(\text{since } (1,3), (2,3), (2,4), (3,4)\in S ).
\end{align*}
The Hasse diagram for $P_{\hfam, S}^{5}$ is as follows.
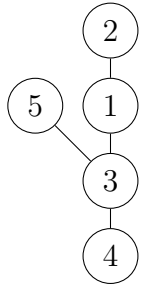
\begin{center}
    \begin{tikzpicture}
        \draw (1,0) node[draw,circle] (v1) {$4$};
        \draw (1,1) node[draw,circle] (v2) {$3$};
        \draw (0,2) node[draw,circle] (v3) {$5$};
        \draw (1,2) node[draw,circle] (v4) {$1$};
        \draw (1,3) node[draw,circle] (v5) {$2$};

        \draw (v1)--(v2);
        \draw (v2)--(v3);
        \draw (v2)--(v4);
        \draw (v4)--(v5);
    \end{tikzpicture}
\end{center}
We get that
\[
\Hinv(S,\hfam(m)) = \Hinv(S,5) = \{34215, 35214, 45213\}.
\]
The inverses of these permutations are precisely the linear extensions of $P_{\hfam, S}^5$,
\[
\Ext(P_{\hfam, S}^5) = \{ 43512, 43152, 43125 \}.
\]
We have $b_k(S)=1=h_{k-1}(P_{\hfam, S}^{5}, 5)$ for each of $k=3,4,5$ and $b_k(S)=0=h_{k-1}(P_{\hfam, S}^{5}, 5)$ otherwise.
\end{example}

As an immediate consequence of Stanley's Theorem~\ref{thm:heightpoly_logconcave} and Lemmas~\ref{lem:height_internal_zeros} and~\ref{lem:Bset_poset_map}, we obtain the main result:

\begin{proposition}\label{prop:b_logconcave} Let $S$ be an $\hfam$-admissible set. The sequence of non-negative integers   $\left(b_k(S)\right)_{k=\hfam(m)-m}^{\hfam(m)}$  is log-concave and does not have internal zeros.  Consequently, the same is true for the sequence $\left(b_{\hfam(m)-\ell}(S)\right)_{\ell=0}^m$.
\end{proposition}

\begin{remark} \label{rem.GG.proof} Let $b_k(S,n):= \#\{\pi\in \Hinv(S,n) \mid \pi_n = k\}$, so $b_k(S,n)$ counts the number of permutation with $\hfam$-inversions equal to $S$ and last entry equal to $k$. Note that $b_k(S;\hfam(m))=b_k(S)$. When $\hfam= (2,3,4,5,\ldots)$, that is, in the setting of descent polynomials, Gaetz and Gao \cite{Gaetz-Gao2021} give an alternative proof of Proposition~\ref{prop:b_logconcave}. They prove, in general, that the sequence $\left(b_k(S,n) \right)_{k=1,\ldots, n}$ is log-concave by induction on $n$. As they note, their methods avoid use of Stanley's results on linear extensions. Unfortunately the same induction techniques do not immediately apply for more general sequences $\hfam$, as there are more conditions defining the sets $\Hinv(S,n)$. 
\end{remark}

\subsection{Log-concavity for \texorpdfstring{$\left(a_k(S)\right)$}{ak(S)}}

Combining Proposition \ref{prop:b_logconcave}, Lemma \ref{lem:ab_conversion}, and results of Brenti~\cite{Brenti1989}, we can achieve an analogous log-concavity result for the sequence of $a_k(S)$-coefficients from Theorem \ref{thm.a-expansion}.

\begin{remark}\label{rem:PF2}
In \cite{Brenti1989}, the property of a sequence $a = (a_0, a_1, \ldots, a_\ell)$ of non-negative integers being log-concave with no internal zeros is called being a $\mathrm{PF}_2$-sequence. As defined in that work, $\left(a_i\right)_{i=0}^\ell$ is $\mathrm{PF}_2$ if every $2 \times 2$ minor of the matrix 
\[
\left[a_{j-i}\right]_{i,j \in \N} = \begin{bmatrix} 
 a_0 & a_1 & a_2 & \cdots & a_n & 0 & \cdots \\
 0 & a_0 & a_1   & \cdots & a_{n-1} & a_n & \cdots \\
 0 & 0 & a_0     & \cdots & a_{n-2} & a_{n-1} & \cdots \\
 \vdots &\vdots&  \vdots  & \ddots  &&
\end{bmatrix}
\] is non-negative, where we set $a_{j-i} = 0$ when $j-i \notin \{0,1,\ldots,n\}$. In particular, this states $\det\left[\begin{matrix}
    a_j & a_{j+k} \\
    a_i & a_{i+k}
\end{matrix}\right] \geq 0$ for all $i<j$ and $k \in \N$. If we let $j=i+1$ and $k=1$, this determinant gives precisely the log-concavity condition. The fact that these conditions are equivalent is proved in \cite[Prop 2.5.1]{Brenti1989}. 
\end{remark}
\begin{proposition}\label{prop:a_logconcave} Let $\hfam$ be a Hessenberg sequence, $S$ be an $\hfam$-admissible set, and $m=\bm(S)$. The sequence of non-negative integers  $(a_i(S))_{i=0}^{m}$ is log-concave and does not have internal zeros.
\end{proposition}
\begin{proof} Since $S$ is fixed, we write $a_k:=a_k(S)$ and $b_k := b_{k}(S)$ throughout the proof. By Lemma \ref{lem:ab_conversion}, we have that 
\[
a_k = \sum_{j=k}^m {j-1 \choose k-1} b_{\hfam(m)-j}
\]
for all $1\leq k \leq m$. By Proposition~\ref{prop:b_logconcave}, the sequence $\left(b_{\hfam(m)-\ell}\right)_{\ell=0}^m$ is log-concave, so by \cite[Thm. 2.5.3]{Brenti1989} the sequence $\left( a_k\right)_{k=1}^{m}$ is log-concave and has no internal zeros.

To complete the proof, we need only to check specifically that $a_1^2 \geq a_0a_2$. To ease notation for the remainder of the proof, let $c_\ell \coloneqq b_{\hfam(m) - \ell}$. By Lemma \ref{lem:ab_conversion}, we have that $a_0 = c_0$, that 
\[
a_1 = \sum_{j=1}^m {j-1 \choose 0} c_{j} = \sum_{j=1}^m c_{j},
\]
and that 
\[
a_2 = \sum_{j=2}^m {j-1 \choose 1} c_j = \sum_{j=2}^m (j-1) c_j.
\]
Since the sequence $\left( c_\ell\right)_{\ell=0}^m$ is log-concave and has no internal zeros, by \cite[Prop. 2.5.1]{Brenti1989} we have that $\det\left[ \begin{matrix}
     c_{j-i} & c_j \\ c_0 & c_i 
\end{matrix}\right] \geq 0$, so $c_0c_j \leq c_{j-i}c_i$ for all $1\leq i < j$.  Thus,
\begin{align*}
        a_0a_2 &= \sum_{j=2}^m (j-1) c_0c_j
         \leq \sum_{j=2}^m \sum_{i=1}^{j-1}c_{j-i}c_i
         \leq \sum_{i,j \in [m]}c_{i}c_j
        = a_1^2.
\end{align*}
This proves $\left(a_k\right)_{k=0}^m$ is a log-concave sequence with no internal zeros. 
\end{proof}

\begin{example}
If $\hfam = (5,5,6,6,7,7,7,\ldots)$ and \[S = \{(1,3),\,(1,4),\,(1, 5),\, (2, 3),\, (2, 4),\, (2, 5),\, (3, 4)\},\] then $m=\bm(S)=3$, $\hfam(m) = 6$, and so 
\begin{align*}
        \hinv(S;n) &=  {n-3 \choose 6-3} + {n-4 \choose 6-4}+{n-5 \choose 6-5} \\
        &= 12{n-6+1\choose 0} + 15 {n-6+1\choose 1} +  6{n-6+1\choose 2}.
\end{align*}
Both the sequences $(1,1,1)$ of $b_k$-coefficients and $(12,15,6)$ of $a_k$-coefficients are log-concave. Note, however, that neither of these sequences satisfy the stronger condition of being the coefficient sequences of a real-rooted polynomial.
\end{example}

\section{A \texorpdfstring{$q$}{q}-analogue of the \texorpdfstring{$\hfam$}{h}-inversion polynomial}

In this section we define a $q$-analogue of $\hinv(S;n)$ and prove $q$-analogues of the expansions in Corollaries~\ref{cor.polynomial} and~\ref{cor.k-expansion}. This section is motivated by work of Gaetz and Gao \cite{Gaetz-Gao2021}, who define a $q$-generalization of the descent polynomial that tracks the length of each permutation. 

Recall that the \emph{length} of a permutation $\sigma\in S_n$ is $\ell(\sigma) = \# \inv(\sigma)$.
Generalizing the definition for $q$-descent polynomials appearing in~\cite{Gaetz-Gao2021}, it is natural to define the \emph{graded $\hfam$-inversion polynomial} as 
\[
\hinv(S, n; q) \coloneqq \sum_{w \in \Hinv(S,n)}q^{\ell(w)}.
\]
By definition, $\hinv(S,n;1) = \hinv(S;n)$. Unlike the polynomial $\hinv(S;n)$, the expression $\hinv(S, n; q)$ is a polynomial in $q$.  The main results of this section, Theorem~\ref{thm:t_graded} and Theorem~\ref{thm:b_graded}, compute formulas for $\hinv(S, n; q)$ which are, respectively, valid for all $n\geq \bj(S)$ and valid when  $\hfam$ is a Hessenberg sequence and $n\geq \hfam(m)$.

\begin{example} Suppose $\hfam(i) = i+2$ for all $i$ and $S = \{(1,3),(2,3),(2,4),(3,4)\}$ as in Example~\ref{ex.poset}, where we computed that 
\[ 
 \Hinv(S,5) = \{34215, 35214, 45213\}.
\]
In this case we have $\hinv(S,5;q) = q^5+q^6+q^7$.
\end{example}

To generalize our expansions for $\hinv(S;n)$ to the graded case, we utilize the \emph{$q$-binomial coefficients}: 
\[
\qchoos{n}{k} \coloneqq \frac{[n]!_q}{[n-k]!_q[k]!_q}\;\;\;\text{ where }\;\;\; [n]!_q = (1)(1+q)\cdots(1+q+\cdots+q^{n-1}).
\]
For example, 
\[
\qchoos{5}{2} 
= \frac{(1+q+q^2+q^3)(1+q+q^2+q^3+q^4)}{(1+q)} = (1+q^2)(1+q+q^2+q^3+q^4).
\]
There are many combinatorial interpretations for $q$-binomial coefficients. We will use the following one (see \cite[\S 1.7]{ec1e2}). 
Define the \emph{length of a subset} $A \subseteq [k,n]$ to be the number of pairs of elements $(a,b)$ in $[k,n]$ such that $a \in A$, $b \notin A$, and $a > b$. Formally,
\[\ell_{[k,n]}(A) \coloneqq \# \left\{ (a ,b) \in [k,n]^2 \mid a>b,\; a \in A,\;b \notin A\right\}.\]

The $q$-binomial coefficient can be written,
\begin{eqnarray} \label{eqn.binom}
\qchoos{n-k+1}{r} = \sum_{\substack{A \subset [k,n]\\\#A = r}} q^{\ell_{[k,n]}(A)}.
\end{eqnarray}
\begin{example}\label{ex:qbinom}
Consider $[2,6] = \{2,3,4,5,6\}$ and subsets $A$ of size $r =2$. The following table shows the lengths (with respect to the superset $[2,6]$) of each such subset.
    \begin{center}
    \begin{tabular}{|c|c|}\hline 
        $\ell_{[2,6]}(A)$ & $A$ \\[1pt]\hline \hline
        $0$ & $\{2,3\}$ \\[2pt]
        $1$ & $\{2,4\}$ \\[2pt]
        $2$ & $\{2,5\}$, $\{3,4\}$  \\[2pt]
        $3$ & $\{2,6\}$, $\{3,5\}$  \\[2pt]
        $4$ & $\{3,6\}$, $\{4,5\}$  \\[2pt]
        $5$ & $\{4,6\}$ \\[2pt]
        $6$ & $\{5,6\}$ \\\hline
    \end{tabular}
    \end{center}
Thus, $\qchoos{5}{2} = 1+q+2q^2+2q^3+2q^4+q^5+q^6$, confirming our earlier computation.
\end{example}

Recall from the proof of Proposition~\ref{prop.polynomial} in Section~\ref{sec.polynomality} that for each $\sigma\in \Hinv(S, \bj(S))$ we define
\[
I_\hfam^{\sigma}(S,n) = \{ \pi \in \Hinv(S,n) \mid \pi|_{\bj(S)} = \sigma  \}.
\]
Recall also that $\bt(\sigma) = \max\{\sigma_r \mid (r,\bj(S)+1)\in \cp_\hfam\}$.

\begin{lemma}\label{lem:tq_expansion} Let $S$ be an $\hfam$-admissible set, $\sigma\in \Hinv(S,\bj(S))$ and $\pi \in I_\hfam^{\sigma}(S,n)$. If $\bt(\sigma)= k$ then  
\[
\ell(\pi) = \ell(\sigma) + \ell_{[k+1,n]}\left(\{\pi_{\bj(S)+1},\ldots,\pi_n\}^c\right)
\]
where $\{\pi_{\bj(S)+1},\ldots,\pi_n\}^c = [k+1, n]\setminus \{\pi_{\bj(S)+1},\ldots,\pi_n\}$.
\end{lemma}
\begin{proof}
This follows from considering how the length of a permutation interacts with the bijection in Equation~\eqref{eqn.polyproofbijection}. First, since $\pi_{\bj(S)+1} < \cdots < \pi_n$ any inversion of $\pi$ is of the form $(i , j)$ with $i \leq \bj(S)$. Any inversion $(i,j)$ of $\pi$ such that $i < j \leq \bj(S)$ is also an inversion of $\pi|_{\bj(S)} = \sigma$.  

%Since $k = \bt(\sigma)$ is the maximum value $\sigma_r$ such that $(r,\bj(S)+1) \in \cp_\hfam$, and so $ k < \pi_{\bj(S)+1}$. 
%In the proof of Proposition~\ref{prop.polynomial} it is shown that $\{\pi_{\bj(S)+1}, \pi_{\bj(S)+2}, \ldots, \pi_n\} \subseteq [k+1, n]$$.
If $(i,j)$ is an inversion of $\pi$ and $i \leq \bj(S) < j$, then $\pi_j$ is an element of the set $\{\pi_{\bj(S)+1},\ldots,\pi_n\}\subseteq [k+1, n]$. In particular, $\pi_i \in [k+1,n]$ as well since $\pi_i > \pi_j$. Thus any inversion $(i,j)$ of $\pi$ such that $i\leq \bj(S)<j$ corresponds to a unique pair $(\pi_i , \pi_j) \in [k+1, n]\times [k+1, n]$ such that $\pi_i>\pi_j$, $\pi_i \in \{\pi_{\bj(S)+1},\ldots,\pi_n\}^c$, and  $\pi_j \in \{\pi_{\bj(S)+1},\ldots,\pi_n\}$. Similarly, any such pair determines a unique inversion.  It follows that the number of these inversions is precisely $\ell_{[k+1,n]}\left(\{\pi_{\bj(S)+1},\ldots,\pi_n\}^c\right)$.
\end{proof}

The following is a direct generalization of~\cite[Theorem 6]{Gaetz-Gao2021} (see Remark~\ref{rem.descentexp} above) and gives a $q$-analogue of the expansion for $\hinv(S;n)$ from Corollary~\ref{cor.polynomial}.

\begin{theorem}\label{thm:t_graded}
Suppose $\hfam$ is a sequence of positive integers such that $\hfam(i)>i$ for all $i$ and let $S$ be an $\hfam$-admissible set. For 
each $n\geq \bj(S)$, we have
\begin{eqnarray}\label{eqn.t_graded_expansion}
\hinv(S,n;q) = \sum_{k=1}^{\bj(S)} t_k(S;q) \qchoos{n-k}{ \bj(S)-k}
\end{eqnarray}
where $\displaystyle t_k(S;q) \coloneqq \sum_{\sigma \in T_k(S,\,\bj(S))} q^{\ell(\sigma)}$ and $T_k(S,\,\bj(S)) = \{\sigma\in \Hinv(S, \bj(S)) \mid \sigma_{\bt(S)} = \bt(\sigma) = k\}$ as in Corollary~\ref{cor.polynomial}.
\end{theorem}
\begin{proof} As in the proof of Proposition~\ref{prop.polynomial}, we consider the decomposition $\Hinv(S,n) = \bigsqcup \Hinv^{\sigma}(S, \bj(S))$ where the disjoint union is taken over all $\sigma\in \Hinv(S,\bj(S))$. Here $\Hinv^{\sigma}(S, \bj(S))$ is the fiber of the flattening map $\Hinv(S,n)\to \Hinv(S, \bj(S))$ over $\sigma$. By the proof of Proposition~\ref{prop.polynomial}, each such fiber is in bijection with subsets of $[\bt(\sigma)+1,n]$. Lemma~\ref{lem:tq_expansion} describes how the bijections in Proposition~\ref{prop.polynomial} interact with the length statistic on the permutations in $\Hinv^{\sigma}(S,n)$. Thus,
\begin{align*}
\hinv(S,n;q) &= \sum_{\sigma\in \Hinv(S,\, \bj(S))} \left(\sum_{\pi \in \Hinv^\sigma(S,n)} q^{\ell(\pi)}\right) \\
&= \sum_{\sigma\in \Hinv(S,\, \bj(S))} \left[ q^{\ell(\sigma)} \left( \sum_{\substack{A \subseteq [\bt(\sigma)+1,n]\\\#A = n-\bj(S)}} q^{\ell_{[\bt(\sigma)+1,n]}(A^c)} \right) \right]\\
&=\sum_{\sigma\in \Hinv(S,\, \bj(S))} q^{\ell(\sigma)} \qchoos{n-\bt(\sigma)}{(n-\bt(\sigma))-(n-\bj(S))}  \qquad  \text{ (by equation~\eqref{eqn.binom})}  \\ 
&=\sum_{\sigma\in \Hinv(S,\, \bj(S))} q^{\ell(\sigma)} \qchoos{n-\bt(\sigma)}{\bj(S) - \bt(\sigma)}.   
\end{align*}
Using the set $T_k(S,\bj(S))$ we can rewrite the above as
\[
\hinv(S,n;q) = \sum_{k=1}^{\bj(S)}  \left( \sum_{\sigma\in T_k(S,\,\bj(S))} q^{\ell(\sigma)}\right) \qchoos{n-k}{\bj(S) - k}.
\]
This completes the proof.
\end{proof}

\begin{example}\label{ex:tbq_example}
Let $\hfam(i)=i+3$ for all $i$ and  $S= \{(3,4), (3,5), (3,6), (4,6), (5,6)\}$ as in Example~\ref{ex.exp2}. We have $\bj(S)=6$ and $\bt(S) = 5$ and in this case,
\[
\Hinv(S, 6) =  \{ 126453,  136452, 236451\}.
\]
So $\sigma_{\bt(S)} = 5$ for all $\sigma \in \Hinv(S, 6)$ and we have
\[
t_5(S;q) = q^5 + q^6 + q^7.
\]
We conclude
\begin{align*}
\hinv(S,n;q) &= (q^5+q^6+q^7) \qchoos{n-5}{6-5} \\
&= (q^5+q^6+q^7) \left(1 + \cdots + q^{n-6}\right) \\
&= q^5+2q^6+3q^7 + \cdots + 3q^{n-1}+2q^n+q^{n+1}.
\end{align*}
\end{example}
\begin{example}
The expansion from Theorem~\ref{thm:t_graded} does not require $\hfam$ to be weakly increasing. For example, consider $\hfam = (7,3,4,5,7,8,9,10,\ldots)$ and 
\[
S = \{(1,4),(1,5),(1,6),(2,3),(3,4)\}
\]
as in Example \ref{ex:ab_negative}, where we saw that $\bj(S)=6$ and $\Hinv(S, 6) = \{465123\}$. We get that $t_4(S;q)=q^{10}$, and thus 
\[
\hinv(S,n;q) = q^{10}\qchoos{n-4}{2}.
\]
\end{example}

When $\hfam$ is a Hessenberg sequence, we can generalize the expansion from Corollary~\ref{cor.k-expansion} of $\hinv(S;n)$ to an expansion of $\hinv(S,n;q)$.
Recall from Section~\ref{sec.expansion1} that 
\[
B_k (S, n) = \{\pi\in \Hinv(S,n) \mid  \pi_{\hfam(m)} = k \}.
\]
By Lemma~\ref{lemma.k-decomp}, $B_k(S;n)\neq \varnothing$ only if $\hfam(m)-m\leq k \leq \hfam(m)$. Recall also that, since $\hfam(m)\geq m$ and $\pi_{m+1}<\pi_{m+2}<\cdots< \pi_n$, we know $\{\pi_{\hfam(m)+1}, \pi_{\hfam(m)+2}, \ldots, \pi_n\} \subseteq [k+1, n]$. 

\begin{lemma}\label{lem:bq_expansion} Let $\hfam$ be a Hessenberg sequence, $S$ be an $\hfam$-admissible set, $m=\bm(S)$, and $k$ such that $\hfam(m)-m\leq k \leq \hfam(m)$. If $\pi \in B_k(S,n)$, then 
\[
\ell(\pi) = \ell(\pi|_{\hfam(m)}) + \ell_{[k+1,n]}\left(\{\pi_{\hfam(m)+1},\ldots,\pi_n\}^c\right)
\]
where $\{\pi_{\hfam(m)+1},\ldots,\pi_n\}^c = [k+1, n]\setminus \{\pi_{\hfam(m)+1},\ldots,\pi_n\}$.
\end{lemma}
\begin{proof} First, since $\pi_{\hfam(m)+1} < \cdots < \pi_n$ any inversion of $\pi$ is of the form $(i , j)$ with $i \leq \hfam(m)$. Any inversion $(i,j)$ of $\pi$ such that $i < j \leq \hfam(m)$ is also an inversion of $\pi|_{\hfam(m)}$.  

If $(i,j)$ is an inversion of $\pi$ and $i \leq \hfam(m) < j$, then $\pi_j$ is an element of the set $\{\pi_{\hfam(m)+1},\ldots,\pi_n\}\subseteq [k+1, n]$. In particular, $\pi_i \in [k+1,n]$ as well since $\pi_i > \pi_j$. Thus any inversion $(i,j)$ of $\pi$ such that $i\leq \hfam(m)<j$ corresponds to a unique pair $(\pi_i , \pi_j) \in [k+1, n]\times [k+1, n]$ such that $\pi_i>\pi_j$, $\pi_i \in \{\pi_{\hfam(m)+1},\ldots,\pi_n\}^c$, and  $\pi_j \in \{\pi_{\hfam(m)+1},\ldots,\pi_n\}$. Similarly, any such pair determines a unique inversion.  It follows that the number of these inversions is precisely $\ell_{[k+1,n]}\left(\{\pi_{\hfam(m)+1},\ldots,\pi_n\}^c\right)$.
\end{proof}

Since Corollary~\ref{cor.k-expansion} is another generalization of an expansion of the descent polynomial, the following is another generalization of~\cite[Theorem 6]{Gaetz-Gao2021}.

\begin{theorem}\label{thm:b_graded}
Let $\hfam$ be a Hessenberg sequence, $S$ be an $\hfam$-admissible set, and $m =\bm(S)$. For all $n\geq \hfam(m)$,
\begin{align*}
\hinv(S,n;q) &= \sum_{k=\hfam(m)-m}^{\hfam(m)} b_k(S;q) \qchoos{n-k}{\hfam(m)-k} \\
&= \sum_{\ell =0}^m b_{\hfam(m)-\ell}(S;q) \qchoos{n-\hfam(m)+\ell }{\ell}
\end{align*}
where $\displaystyle b_k(S;q):=  \sum_{\pi \in B_k(S,\hfam(m))} q^{\ell(\pi)}$.
\end{theorem}
\begin{proof} Lemma~\ref{lem:bq_expansion} describes how the bijection from the statement of Lemma~\ref{lemma.k-decomp} interacts with the length statistic on the permutations in $B_k(S,n)$. Thus,
\begin{align*}
\hinv(S,n;q) &= \sum_{k=\hfam(m)-m}^{\hfam(m)} \left(\sum_{\pi \in B_k(S,n)} q^{\ell(w)}\right) \\
&= \sum_{k=\hfam(m)-m}^{\hfam(m)} \left[ \left( \sum_{\pi \in B_k(S,\hfam(m))} q^{\ell(\pi)} \right) \left( \sum_{\substack{A \subseteq [k+1,n]\\\#A = n-\hfam(m)}} q^{\ell_{[k+1,n]}(A^c)} \right) \right]\\
&= \sum_{k=\hfam(m)-m}^{\hfam(m)} b_k(S;q) \qchoos{n-k}{(n-k)-(n-\hfam(m))}  \qquad  \text{ (by equation~\eqref{eqn.binom})}  \\
&= \sum_{k=\hfam(m)-m}^{\hfam(m)} b_k(S;q) \qchoos{n-k}{\hfam(m)-k}.  
\end{align*}
%The final equality in the theorem follows by shifting the indexing of the sum, from $k$ to $\ell = \hfam(m)-k$.
This completes the proof.
\end{proof}

\begin{example}\label{ex:graded_h_inversion_poly} Let $\hfam(i)=i+3$ for all $i$ and  $S= \{(3,4), (3,5), (3,6), (4,6), (5,6)\}$ as in Example~\ref{ex.exp2} and Example~\ref{ex:tbq_example}. We have $m=\bm(S)=5$ and $\hfam(m)=8$ and in this case,
\begin{align*}
\Hinv(S, 8) =  \{ 12645378, &12745368, 12845367, 13645278, 13745268, \\& \quad\quad\quad13845267, 23645178, 23745168, 23845167\}.
\end{align*}
Therefore we have:
\begin{align*}
B_{7}(S,8) &= \{ 12845367, 13845267,23845167 \}  \\
B_8(S,8) &= \{12645378,12745368, 13645278, 13745268, 23645178,23745168   \}.
\end{align*}
Thus,
\[
b_7(S;q) = q^7+q^8+q^9 \quad \text{ and } \quad b_8(S;q) = q^5+2q^6+2q^7+q^8,
\]
and we conclude
\begin{align*}
\hinv(S,n;q) &= (q^7+q^8+q^9) \qchoos{n-7}{8-7} +  (q^5+2q^6+2q^7+q^8) \qchoos{n-8}{8-8} \\
&= (q^7+q^8+q^9)(1+q+\cdots + q^{n-8}) + q^5+2q^6+2q^7+q^8. 
\end{align*}
The reader can confirm that when $q=1$, we recover the expansion from Example~\ref{ex.exp2}, and that when evaluating at $n \geq 8$ that we get the same expansion as in Example~\ref{ex:tbq_example}. Notice also that this formula is not necessarily valid as written for values of $n<\hfam(m)$, even though we know an expansion exists when $\bj(S) \leq n < \hfam(m)$ because of Theorem \ref{thm:t_graded}.
\end{example}

In \cite[Theorem 6]{Gaetz-Gao2021} Gaetz and Gao show that, in case of graded descent polynomials (namely $\hfam(i)=i+1$), the sequences of polynomials $\left(t_k(S;q)\right)_{k=1}^{\bj(S)}$ and $\left( b_{\hfam(m)-\ell}(S;q) \right)_{\ell=0}^{m}$ from the expansions in Theorem~\ref{thm:t_graded} and Theorem~\ref{thm:b_graded} (which agree for this $\hfam$) are strongly $q$-log concave sequences.  It is an open question whether this fact generalizes to arbitrary graded $\hfam$-inversion polynomials. We have confirmed that this is the case for $\left(t_k(S;q)\right)_{k=1}^{\bj(S)}$ whenever $\bj(S) \leq 6$, and the same for $\left( b_{\hfam(m)-\ell}(S;q) \right)_{\ell=0}^{m}$ when $\hfam$ is weakly increasing and $\hfam(\bm(S))\leq 8$, and so make the following conjecture.

\begin{conjecture}\label{conj:strongly_q_logconcave}
Let $S$ be an $\hfam$-admissible set. The polynomial sequence $(t_k(S;q))_{k=1}^{\bj(S)}$ is strongly $q$-log concave. If $\hfam$ is a Hessenberg sequence and $m =\bm(S)$, then the polynomial sequence $(b_k(S;q))_{k=\hfam(m)-m}^{\hfam(m)}$ is strongly $q$-log concave.
\end{conjecture}
As noted in Remark~\ref{rem.GG.proof} above, the inductive proof techniques from~\cite{Gaetz-Gao2021} do not generalize easily to our setting. Similarly, it is not obvious how to formulate a graded variation of the expansion from Theorem~\ref{thm.a-expansion}. It would also be interesting to characterize which (labeled) finite posets $P$ have a graded height polynomial $ \sum_{\phi \in \text{Ext}(P)} q^{\ell(\phi)}x^{\phi(v)-1}  = \sum_{k=0}^{\#P-1} h_k(P,v;q)\,x^k$ such that the coefficient sequence $(h_k(P,v;q))$ for each $v\in P$ is strongly $q$-log concave. We leave these questions for consideration in future research.

\subsection*{Acknowledgments} The authors are supported in part by NSF CAREER grant DMS-2237057. The first author was also supported by the Freiwald Scholars Program at Washington University in St.~Louis. The authors are grateful to the anonymous referee for insightful remarks which improved the exposition and strengthened the results of this paper.

%\nocite{*}
%\bibliographystyle{acm}
%\bibliography{references.bib}

\end{document}